\documentclass[twoside]{article}
\usepackage{graphicx} % Required for inserting images
\usepackage{geometry}
\usepackage{mathtools}
\usepackage{amsmath}
\usepackage{amssymb}
\usepackage{bm}
\usepackage{xcolor}
\usepackage{caption}
\usepackage{float}
\usepackage{algpseudocode}
\usepackage{comment}
\usepackage{enumitem}
\numberwithin{equation}{section}

\newtheorem{remark}{Remark}[section]

\usepackage{titlesec}

\newcommand{\bfu}{\mathbf{u}}
\newcommand{\bfv}{\mathbf{v}}
\newcommand{\bfs}{\mathbf{s}}
\newcommand{\Mphi}{\mathbf{M}_{\mphi}}
\newcommand{\Kphi}{\mathbf{K}_{\mphi}}

\newcommand{\x}{\mathbf{x}}

\newcommand{\mphi}{\bm{\phi}}

\newcommand{\scalp}[2]{\langle #1\, , #2 \rangle}

\newcommand{\dtheta}{\delta\theta}

\newcommand{\bfa}{\mathbf{a}}

\usepackage{RR}
\RRNo{9583}
\RRdate{March 2025} 
\RRauthor{
Julien Besset\thanks[makutu]{Makutu, Inria, TotalEnergies, University of Pau et des Pays de l'Adour, 64000 Pau, France.}
\and
H\'el\`ene Barucq\thanksref{makutu}
\and
Rabia Djellouli\thanks[iris]{Department of Mathematics and Interdisciplinary Research Institute for the Sciences (IRIS), California State University, Northridge, CA 91330, USA.}
\and
Stefano Frambati\thanksref{makutu}
}

\authorhead{J. Besset \& H. Barucq \& R. Djellouli \& S. Frambati}
\RRtitle{MOR-T$_L$ : Une nouvelle méthode de réduction d'ordre pour les problèmes paramétrés avec application à la propagation des ondes sismiques.}

\RRetitle{{MOR-T$_L$: A Novel Model Order Reduction Method for Parametrized Problems with Application to Seismic Wave Propagation}}
\titlehead{MOR-T$_L$: A Novel Model Order Reduction Method for Parametrized Problems}
\RRresume{Ce rapport présente MOR-T$_L$, une nouvelle méthode de réduction d’ordre basée sur le développement de Taylor et les dérivées de Fréchet par rapport aux paramètres du modèle. Cette approche permet de construire des bases réduites avec un coût computationnel minimal à l'aide d'une décomposition QR tout en assurant une robustesse accrue face aux variations des paramètres du modèle.\\
L’efficacité de la méthode MOR-T$_L$ est démontrée sur des problèmes de propagation d’ondes sismiques en deux dimensions. Nos expériences numériques montrent que cette méthode permet de réduire significativement le coût de calcul lié à la construction de la base réduite par rapport à une approche classique comme la décomposition orthogonale aux valeurs propres, tout en maintenant une précision élevée, y compris pour des variations du paramètre du modèle importantes.\\
Ces résultats positionnent MOR-T$_L$ comme une alternative prometteuse pour les simulations intensives et les problèmes inverses, notamment en inversion sismique, où des mises à jour répétées du modèle sont nécessaires.}
\RRabstract{This paper presents an efficient strategy for constructing Reduced-Order Model (ROM) bases using Taylor polynomial expansions and Fréchet derivatives with respect to model parameters. The proposed approach enables the construction of ROM bases with minimal additional computational cost. By exploiting Fréchet derivatives -- solution to the same problem with distinct right-hand sides -- the method introduces a streamlined multiple-right-hand-side (RHS) strategy for ROM bases construction. This approach not only reduces overall computational expenses but also improves accuracy during model parameter updates. Numerical experiments on a two-dimensional wave problem demonstrate significant efficiency gains and enhanced performance, highlighting the potential of the proposed method to advance computational cost-effectiveness, particularly in seismic inversion applications. 
}

\RRmotcle{Model Order Reduction (MOR) -- Reduced Order Model (ROM) -- Fréchet derivatives
 -- Taylor polynomial representation -- Spectral Element Method (SEM) -- QR decomposition
 -- Line Search}
\RRkeyword{
 Méthode de Réduction d'Ordre (MOR) -- Modèle d'Ordre Réduit (ROM) -- Dérivées de Fréchet
 -- Polynôme de Taylor -- Elements Finis Spectraux (SEM) -- Décomposition QR
 -- Recherche Linéaire
}

\RRprojet{Makutu}  % cas d'un seul projet
\RCBordeaux 

\begin{document}
\RRNo{9583} 
\makeRR    % cas d'un rapport de recherche
%% \makeRT % cas d'un rapport technique.
%% a partir d'ici, chacun fait comme il le souhaite

\tableofcontents 
\newpage
% ------        main sections ------------------------------------------

\section{Introduction and Motivation}
The development of efficient solution methodologies for mathematical models that describe complex, real-world phenomena is of paramount importance across a wide range of fields, including heat transfer in materials with varying thermal properties, pollutant transfer in fluids, semiconductor applications, seismic and geophysical problems, radar and sonar technologies, and climate dynamics. These phenomena are typically governed by Partial Differential Equations (PDEs) whose discretization results in large, high-dimensional systems that are computationally expensive to solve. Moreover, these models are often highly sensitive to variations in multiple parameters, such as material properties, environmental conditions, or geometrical factors, and there is frequently a need to solve these problems repeatedly for different parameter values. This is especially true in inverse problems, where one seeks to reconstruct a set of parameters (or all of them) from experimental measurements, requiring the repeated solution of direct problems. As the number of parameters increases, the computational burden becomes prohibitive, making this process practically unattainable.\\
Over the past fifty years, significant efforts have been dedicated to the design of numerical strategies that strike a balance between computational cost and accuracy. Approaches such as adaptive grids \cite{diez1999unified, zienkiewicz1991adaptivity}, where mesh sizes are adjusted based on parameter variations, and variable approximation orders \cite{dumbser2006arbitrary, beriot2019anisotropic}, where the order of approximation changes in response to parameter changes, have been employed to reduce the computational load. Additionally, methods like Trefftz \cite{herrera2000trefftz, despres2022trefftz} and hybridizable discontinuous Galerkin (HDG) methods \cite{cockburn2012discontinuous, cockburn2000development, di2011mathematical, hesthaven2007nodal} have been introduced, which use fewer degrees of freedom by exploiting local solutions or static condensation techniques. Although these methods offer computational advantages, their performance remains unsatisfactory in situations where model parameters change frequently, especially in the context of model updates where accuracy must be preserved while also avoiding dramatic increases in computational cost.\\
A promising recent development in this domain is the use of Model Order Reduction (MOR) techniques, which aim to approximate the high-dimensional solution space of a system using a significantly smaller subspace, thereby reducing the size of the system and improving computational efficiency. Well-known MOR techniques include Proper Orthogonal Decomposition (POD) \cite{cordier_reduction_2006}, Dynamic Mode Decomposition (DMD) \cite{tu_dynamic_2014}, and neural network-based ROM \cite{mohan_deep_2018}. However, the major challenge with these approaches is their sensitivity to parameter changes. As parameters vary, the accuracy of the ROM-based solutions can degrade significantly. The simplest ways to address this issue involve reconstructing or enlarging the basis functions, both of which result in increased computational costs and undermine the very efficiency gains that ROM aims to achieve. Alternative methods, such as Proper Generalized Decomposition (PGD) \cite{chinesta_recent_2010} and Grassmann Interpolation \cite{meza_interpolation_2018}, have been proposed, but their improvements in accuracy and computational performance have been limited.\\
In response to these challenges, we propose a novel MOR approach based on Taylor polynomial expansions. The key innovation of this method lies in how the ROM basis functions are constructed: we use Taylor polynomial representations to generate these bases, which leads to a more robust method that can handle parameter variations more effectively. We have named this new approach {\bf MOR-T$_L$} where T stands for Taylor polynomials and L represents the polynomial degree. The proposed MOR-T$_L$ strategy provides a unique combination of advantages: it reduces the size of the system, thereby lowering computational cost, while also being robust to parameter changes. This dual benefit is critical not only for maintaining accuracy when the model is updated but also for preventing the computational cost from rising dramatically as parameters vary.\\
The theoretical foundation and numerical implementation of MOR-T$_L$ are presented in this paper. We demonstrate its efficiency in reducing computational time while maintaining high accuracy, even as parameters change. By leveraging the Fréchet derivatives of the solution with respect to the problem parameters - solutions to the same problem with distinct right-hand sides - we enable a streamlined multiple-right-hand-side (RHS) strategy for constructing the ROM bases. This method allows the adaptive construction of reduced bases with minimal additional computational cost. As we show in Section 5, MOR-T$_L$ is highly effective in addressing parameter changes, making it a promising candidate for solving large-scale problems in a variety of applications.\\
The remainder of this paper is organized as follows. In Section 2, we introduce the adopted notations, assumptions, and the class of parametrized partial differential equations under consideration. Section 3 presents the theoretical foundation of the MOR-T$_L$ method and discusses its computational complexities and cost. In Section 4, we apply MOR-T$_L$ to a parametrized two-dimensional wave problem arising in seismic applications. Finally, Section 5 provides a numerical performance comparison of MOR-T$_L$ with the standard POD approach and the Spectral Element Method (SEM), highlighting the significant computational savings and robustness of MOR-T$_L$, which make it a strong candidate for large-scale seismic applications.

\section{Preliminaries: Notations, Assumptions, and Problem Setup}
In this section, we present the adopted notations, assumptions, and the class of boundary value problems under consideration.

\subsection{Nomenclature and Assumptions}
Throughout this paper, we adopt the following notations and hypothesis:
\begin{itemize}
    \item $\x$ represents the spatial variable. $\x=(x_1,x_2) \in \mathbb{R}^2$.
    \item $t$ denotes a time variable. $t \geq 0$.
    \item $N_x$ is a natural number representing the number of degrees of freedom associated with the discretization scheme.
    \item $N_t$ is a natural number representing the number of time steps.
    \item $\Delta t$ is a positive number representing the time step.
    \item $\x_j$, $j=1,\ldots,N_x$.
    \item $t_j=j\Delta t$, $j=1,\ldots, N_t$.
    \item $T$ is the timelength of the computations. $T= N_t \Delta t$.
    \item $\Omega$ is a bounded domain in $\mathbb{R}^2$ whose boundary $\partial \Omega$ is Lipschitz continuous \cite{evans1992studies}.
    \item $X,Y$, and $Z$ are respectively Banach spaces, corresponding to the space of solutions, source terms and model parameters, respectively.
    \item $\theta$ represents the model parameter. $\theta\in Z$.
    \item $f$ is a source term. $f \in Y$.
    \item $u_\theta$ is the solution of the considered boundary value problem with a parameter $\theta$. $u_\theta \in X$.
    \item $\dtheta$ represents the model parameter perturbation direction. $\dtheta\in Z$.
    \item The derivative of $u_\theta$ at $\theta = \theta_0$ in the direction of $\dtheta_0$ is given by:
    \begin{equation}\label{diff}
        Du_{\theta_0}\dtheta_0 = \lim_{\epsilon \rightarrow 0} \frac{1}{\epsilon}\left( u_{\theta_0 + \epsilon\dtheta_0} - u_{\theta_0} \right)
    \end{equation}
    Hence, we denote the Fréchet derivative of order $\ell$ of $u_\theta$ at $\theta = \theta_0$ in the direction of $\dtheta_0$ by $v^\ell_{\theta_0} = D^\ell u_{\theta_0} \dtheta_0^\ell$, and $v^0_{\theta_0} = u_{\theta_0}$.
    \item The Taylor polynomial of degree $L$ of $u_\theta$ centered at $\theta_0$ with a perturbation $\dtheta_0$ is given by:
    \begin{equation}\label{taylor}
        P_{\theta_0}(u_\theta) = \sum_{\ell=0}^L \frac{1}{\ell!} v^\ell_{\theta_0}
    \end{equation}
    \item $\mathcal{L}_\theta$ is a linear operator from $X$ to $Y$. $\mathcal{L}_\theta$ depends on the parameter $\theta$.
    \item The derivative of order $\ell$ of $\mathcal{L}_\theta$ at $\theta = \theta_0$ in the direction of $\dtheta$ is denoted by $\displaystyle D^\ell \mathcal{L}_{\theta_0} \dtheta_0^\ell$, and $D^0\mathcal{L}_{\theta_0}\dtheta_0^0 = \mathcal{L}_{\theta_0}$.
    \item $\delta_{ij}$ defines the Kronecker symbol, i.e. $\delta_{ij}=1$ if $i=j$, and $0$ otherwise.
    \item $\delta(\x-\x_0)$ is the delta Dirac function at $\x=\x_0$.  
\end{itemize}

\subsection{Mathematical Model}
Given $f \in Y$ and for $\theta \in Z$, we consider the following class of linear PDE systems:
%We consider the following class of boundary value problems:

\begin{equation}\label{PDE}
    \begin{cases}
        \text{Find } u_\theta \in X \text{ such that }\\
        \mathcal{L}_\theta u_\theta = f 
    \end{cases}
\end{equation}

The goal is to develop an efficient computational procedure to solve \eqref{PDE} for a range of values of the parameter $\theta$ in a bounded domain. This is particularly important for applications like inverse problems, including inverse seismic scenarios, where solving $\eqref{PDE}$ for varying $\theta$ is often required.\\
It is essential to note that classical representations of the solution $u_{\theta_0}$ using POD for instance, only applies for the fixed value of $\theta = \theta_0$. This limitation arises because the basis functions depend on $\theta_0$, and thus this representation cannot be used to approximate accurately solution to \eqref{PDE} for $\theta \neq \theta_0$. As noted in \cite{jb}, even a small perturbation in $\theta_0$ can lead to significant inaccuracies. Therefore, standard MOR approaches at this stage are inadequate for solving \eqref{PDE} for multiple values of $\theta$. To overcome this challenge, based on Taylor polynomial expansion, we propose to construct a basis for the Fréchet derivatives of $u_\theta$ evaluated at $\theta=\theta_0$ in the direction of perturbation $\dtheta=\dtheta_0$, then reconstruct the solution $u_\theta$ evaluated in $\theta_\alpha = \theta_0 + \alpha\dtheta_0$ with $\alpha\in\mathbb{R}$.

\section{Introducing MOR-T$_L$: A Novel Reduced Order Method for Solving \eqref{PDE} Across Parameter Variations}

The method we propose for solving \eqref{PDE} for a set of admissible parameters $\theta$ falls within the category of MOR methods. What sets our approach apart is the novel way in which we construct the basis functions. Specifically, we use Taylor polynomial representations to build this basis. For this reason, we have named the proposed method MOR-T$_L$: T stands for Taylor polynomials, and $L$ represents the polynomial degree. 
The strength of the MOR-T$_L$ strategy lies in its ability to solve \eqref{PDE} for different values of $\theta$ without sacrificing accuracy. In contrast to existing ROM methods, where the constructed basis for a specific $\theta$ is not suitable for solving \eqref{PDE} when $\theta_\alpha = \theta_0 + \alpha \dtheta_0$ (even for small perturbations $\alpha\dtheta_0$) our approach retains its accuracy. In traditional MOR methods, such perturbations often result in a significant loss of accuracy \cite{jb}. To correct this, the simplest options are to either reconstruct the basis or enlarge it, both of which significantly increase computational costs and undermine the efficiency gains that MOR methods aim to provide.\\

MOR-T$_L$’s approach is to construct the basis in two distinct stages, outlined below:

\subsection*{\textbf{Stage 1: Construction of the Bases for the Fréchet Derivatives Solutions}}

The first step is to reformulate \eqref{PDE}, evaluated in $\theta_\alpha = \theta_0 + \alpha\dtheta_0$, as a series of $L+1$ PDEs using the Taylor polynomial of order $L$ centered at $\theta_0$ given by \eqref{taylor} for $\alpha\in\mathbb{R}$. In the following, we consider the classical Taylor expansion of order $L$ as described by the polynomial defined in \eqref{taylor}.
\begin{equation}\label{PDE taylor}
    \begin{cases}
        \displaystyle \text{Find } v_{\theta_0}^\ell \in X, \; \ell = 0,\ldots,L, \text{ such that }\\
        \displaystyle \mathcal{L}_{\theta_0} \left( v_{\theta_0}^\ell \right) = f^{(\ell)}
    \end{cases}
\end{equation}
where $\displaystyle v^\ell_{\theta_0} = D u_{\theta_0} \dtheta_0$, and $\displaystyle f^{(\ell)} = - \sum_{k=0}^{\ell - 1} {\ell \choose k} \left( D^{\ell - k} \mathcal{L}_{\theta_0}\dtheta_0^{\ell-k} \right)( v_{\theta_0}^k )$ with $\displaystyle {\ell \choose k} = \frac{\ell!}{k!(\ell-k)!}$; $0\leq k \leq \ell-1$.\\ 
Note that $\alpha$ drops in \eqref{PDE taylor} by linearity of $D$ and $\mathcal{L}_{\theta_0}$:
\begin{equation}
    (D^{\ell-k} \mathcal{L}_{\theta_0}(\alpha \dtheta_0)^{\ell-k})(\alpha^k v^k_{\theta_0}) = \alpha^\ell (D^{\ell-k} \mathcal{L}_{\theta_0} \dtheta^{\ell-k})(v^k_{\theta_0})
\end{equation} 
This property is of great practical interest since it allows the reconstruction of the solution of \eqref{PDE} for $\theta = \theta_\alpha$ with different values of $\alpha$ using \eqref{taylor} up to multiplication factors.\\
Also note that $v^\ell_{\theta_0}$ are solutions to the PDE \eqref{PDE}, but with different right-hand sides, resulting in significant computational cost savings.\\

As with any domain-based discretization scheme, if $\eqref{PDE taylor}$ is set in an infinite domain, the problem must first be reformulated in a bounded domain $\Omega$. This is achieved by introducing an artificial external boundary to which a boundary condition is applied. The boundary condition is characterized by the boundary operator $B_\theta$ (see, e.g., \cite{turkel_boundary_2008, givoli_high-order_2004} and the references therein). Therefore, the boundary value problems considered in the bounded domain are given by:

\begin{equation}\label{BVP taylor}
    \textbf{(BVP)}
    \begin{cases}
        \displaystyle \text{Find } v_{\theta_0}^\ell \in X, \; \ell = 0,\ldots,L&, \text{ such that } \\
        \displaystyle \mathcal{L}_{\theta_0} \left( v_{\theta_0}^\ell \right) = f^{(\ell)} &\text{in} \quad \Omega \times (0,T)\\
        \displaystyle \mathcal{B}_{\theta_0} \left( v_{\theta_0}^\ell \right) = 0 &\text{on} \quad \partial\Omega \times (0,T)
    \end{cases}
\end{equation}
Then, the construction of the reduced bases follows a two-steps procedure.

\paragraph{ \it \textbf{(a) Initialization Step.}}
For $\theta = \theta_0$ fixed, we solve \eqref{BVP taylor} using a full order high-fidelity space discretization scheme, e.g., finite element approximation (FEM), that we can represent by the following Galerkin decomposition:

\begin{equation}\label{vl fem}
    \begin{aligned}
        v^\ell_{\theta_0}(\x,t) = \sum_{i=1}^{N_x} v^\ell_{\theta_0,i}(t) \psi_i(\x); \quad \ell=0,\ldots, L
    \end{aligned}
\end{equation}
where $\{\psi_i\}_{i=1}^{N_x}$ represents a nodal finite element basis associated with the considered discretization scheme. To simplify the notation, we omit the explicit reference to the dependence of the approximate solution on the domain discretization scheme. Thus, the notation $v_{\theta_0}^\ell$ should be understood as the solution to the discretized PDE. Let $\bfv^\ell_{\theta_0}$ be the approximate solution vector. Hence:

\begin{equation}\label{vl vector}
    \begin{aligned}
        \bfv^\ell_{\theta_0}(t) = [v^\ell_{\theta_0,1}(t),\ldots,v^\ell_{\theta_0,N_x}(t)]^T; \quad \ell=0,\ldots, L.
    \end{aligned}
\end{equation}
Then, for a sequence of times $t_i$, $i=1,\ldots, N_t$, we define the snapshot vector $\bfs^\ell_i$ as:

\begin{equation}\label{snapshot vector sl}
    \begin{aligned}
        \bfs^\ell_i = \bfv^\ell_{\theta_0}(t_i), \quad i=1,2,\ldots,N_t ; \quad \ell=0,\ldots, L
    \end{aligned}
\end{equation}
which represents the discrete version of the solution field to \eqref{BVP taylor} at time $t_i$. We can then define the snapshot matrix $S^\ell$ whose columns vectors are $\bfs^\ell_i$, viz., 

\begin{equation}\label{Sl}
    S^\ell =
    \begin{aligned}
        \begin{bmatrix}
            \vrule &\vrule & &\vrule \\
            \bfs^\ell_1 &\bfs^\ell_2 &\cdots &\bfs^\ell_{N_t} \\
            \vrule &\vrule & &\vrule
        \end{bmatrix}; \quad \ell=0,\ldots, L
    \end{aligned}
\end{equation}
Note that $S^\ell$ is an $N_x \times N_t$ real matrix. Then, the solution vector $\bfv^\ell_{\theta_0}(t)$ at any given time is represented as a finite sum, which is expressed as:

\begin{equation}\label{basis Sl}
    \begin{aligned}
        \bfv^\ell_{\theta_0}(t) = \sum_{i=1}^{N_t} \Tilde{a}^\ell_i(t) \bfs^\ell_i; \quad \ell=0,\ldots, L
    \end{aligned}
\end{equation}
Here, the time-dependent coefficients $\{\Tilde{a}^\ell_i\}_{i=1}^{N_t}$ are solutions to the system of ordinary differential equations, which emerge when substituting \eqref{basis Sl} into the problem \eqref{BVP taylor}. It is important to note that the field representation in equation \eqref{basis Sl} assumes that $\bfv^\ell_{\theta_0}$ lies in the span of the vectors $\{\bfs^\ell_1,\bfs^\ell_2,\cdots,\bfs^\ell_{N_t}\}$. Hence, the latter need not be linearly independent, so it is crucial to address the potential redundancy within the set and reduce it to retain only the independent vectors. This reduction is necessary to minimize both storage and computational costs. \\
To emphasize the need for this reduction, consider that in seismic applications, the number of time iterations can be very large ($\sim10^4$) and the value of $N_x$ can reach up to $10^9$. Storing all computed snapshots during the simulation becomes computationally prohibitive, if not outright impossible.

\paragraph{\it \textbf{(b) Reduction Step.}}
There are several well-established methods for constructing an orthonormal basis from a given set of vectors, such as Singular Value Decomposition (SVD) \cite{golub_matrix_2013, strang_introduction_2016, trefethen_numerical_1997, eckart_approximation_1936} and QR decomposition methods \cite{golub_matrix_2013, strang_introduction_2016, trefethen_numerical_1997, leon_gramschmidt_2013}. Applying any of these techniques results in the construction of an orthonormal reduced set of modified snapshot vectors which form the reduced basis

\begin{equation}
    \begin{aligned}
        [\mphi^\ell_1, \mphi^\ell_2, \cdots, \mphi^\ell_{M_\ell}]; \quad M_\ell \leq N_t; \quad \ell=0,\ldots, L
    \end{aligned}
\end{equation}
This reduction also leads to the full-rank modified snapshots matrix $\Tilde{S}^\ell$ given by:

\begin{equation}\label{modified Sl}
    \Tilde{S}^\ell = 
    \begin{aligned}
        \begin{bmatrix}
            \vrule &\vrule & &\vrule \\
            \mphi_1^\ell &\mphi_2^\ell &\cdots &\mphi_{M_\ell}^\ell \\
            \vrule &\vrule & &\vrule
        \end{bmatrix}; \quad \ell=0,\ldots, L
    \end{aligned}
\end{equation}
with $M_\ell$ being always less than $N_t$. Hence, we represent $\bfv^\ell_{\theta_0}$ as follows:

\begin{equation}\label{modified basis Sl}
    \begin{aligned}
        \bfv^\ell_{\theta_0}(t) \approx \sum_{i=1}^{M_\ell} \hat{a}^\ell_i(t) \mphi^\ell_i; \quad \ell=0,\ldots, L
    \end{aligned}
\end{equation}

\begin{remark}
Note that MOR-T$_0$ equipped with SVD is the classical POD.
\end{remark}
\subsection*{\textbf{Stage 2: Construction of the reduced model}} \label{stage 3}
This stage consists in selecting the final basis for MOR-T$_L$, following the reduction process and incorporating the Fréchet derivatives.
It follows from \textbf{Stage 1}, the obtention of an augmented modified snapshot matrix given by

\begin{equation}\label{augmented S}
    S^{\mathrm{\mathrm{aug}}} =
    \begin{aligned}
        \begin{bmatrix}
            [\Tilde{S}^0] &[\Tilde{S}^1] &\cdots &[\Tilde{S}^{L}]\\
        \end{bmatrix}
    \end{aligned}
\end{equation}
where the block matrix $\Tilde{S}^0,\ldots,\Tilde{S}^L$ are given by \eqref{modified Sl}.\\
Note that the matrix $S^{\mathrm{\mathrm{aug}}}$ is therefore a $N_x \times \displaystyle \sum_{\ell=0}^{L} M_\ell$ matrix. In addition, each block matrix $\Tilde{S}^\ell$ is a full-rank matrix since the column vectors form an orthonormal basis. However, this does not ensure that $S^{\mathrm{\mathrm{aug}}}$ is a full-rank matrix. Hence, modulo a reordering of the columns of \eqref{augmented S}, we apply again a reduction order procedure and obtain a new orthonormal set of vectors:

\begin{equation}\label{final phi basis}
    \begin{aligned}
        [\mphi_1, \mphi_2, \cdots, \mphi_{N}]; \quad N \leq \sum_{\ell=0}^{L} M_\ell.
    \end{aligned}
\end{equation}
The number $N$ is the result of applying the reduction method to the final augmented snapshot matrix. 
Therefore, for $\theta_\alpha = \theta_0 + \alpha\dtheta_0$, $\bfu_{\theta_\alpha}$ can be approximated, using Taylor polynomial expansions of $u_\theta$ centered at $\theta_0$, as follows:

\begin{equation}\label{modified augmented basis S}
    \begin{aligned}
        \bfu_{\theta_\alpha}(t) &= \sum_{i=1}^{N} a_i(t) \mphi_i
    \end{aligned}
\end{equation}
The computation of $\{a_i\}_{i=1}^N$ requires solving the ODE system resulting from substituting \eqref{modified augmented basis S} into the space-discretized problem \eqref{PDE} for the new value $\theta_\alpha$.
Consequently, the proposed representation \eqref{modified augmented basis S} is used to compute $u_\theta$ for arbitrary values of $\theta$ along the direction of $\dtheta$, in a substantially faster way by solving a reduced system.

\begin{remark}
The following points are worth nothing:\\
%\begin{enumerate}[label=\alph*), leftmargin=*]
{\bf{a.}} The right-hand sides of the boundary value problems in \eqref{BVP taylor} do not depend on the source term $f$ because $f$ is independent of the parameter $\theta$. If $f$ depends on $\theta$, the method still applies, but with updated right-hand sides that include the successive Fréchet derivatives of $f$ with respect to $\theta$ up to order $L-1$.\\
{\bf{b.}} MOR-T$_L$ can also handle perturbations of the parameter in multiple directions of the form: $\theta_\alpha = \theta_0 + \displaystyle \sum_{k=1}^K \alpha_k \dtheta_k$ with $\scalp{\dtheta_i}{\dtheta_j}=0$, $\forall i \neq j$. However, this introduces the need to compute the Fréchet derivatives for the cross-directions, thereby increasing the number of problems to solve combinatorially.
%\end{enumerate}

\end{remark}

\subsection*{\textbf{Final remark on the Computational costs and complexity of the MOR-T$_L$}}
In this section, we provide an overview of the computational costs associated with the MOR-T$_L$ method, which can be divided into two categories: (i) the Offline computational cost, which covers the time required to construct the reduced basis at a fixed parameter value $\theta_0$ as well as the reduced order system itself, and (ii) the Online computational cost, which accounts for the computation of the reduced-order system solution and the reconstruction of the numerical solution for an updated parameter $\theta_\alpha$.

\subsubsection*{Offline Computational Cost}
The offline computational cost primarily involves solving the discrete system with multiple right-hand sides and construct the reduced-order basis. These tasks are considered pre-processing steps since they are performed once for a fixed parameter $\theta_0$. Specifically, the components of the offline cost are as follows:
\begin{itemize}
    \item \textbf{Construction of $L + 1$ snapshots matrices $S^\ell$ \eqref{Sl}:} This task requires solving a full boundary value problem with multiple right hand sides. The computational cost depends on the chosen numerical scheme for solving $\eqref{BVP taylor}$ and the solver (direct of iterative) employed for the resulting system. 
    \item \textbf{Construction of $L + 1$ modified (orthonormal) snapshots matrices $\Tilde{S}^\ell$ \eqref{modified Sl}:} This step involves a reduction technique. When using randomized SVD \cite{golub_matrix_2013, strang_introduction_2016, trefethen_numerical_1997}, the computational cost is of order $\displaystyle \sum_{\ell=0}^L \mathcal{O}(N_x N_t (M_\ell + p) + (N_x+N_t)(M_\ell+p)^2)$ with $p$ an oversampling, whereas using QR decomposition with the Gram-Schmidt orthonormalization procedure \cite{golub_matrix_2013, strang_introduction_2016, trefethen_numerical_1997, leon_gramschmidt_2013} results in a cost of order $\displaystyle \sum_{\ell=0}^L \mathcal{O}(N_x M_\ell)$.
    \item \textbf{Construction of the final basis \eqref{final phi basis}:} This requires applying a reduction technique (SVD or QR) to the augmented modified snapshot matrix \eqref{modified augmented basis S}. The computational cost of this task is similar to the previous step. 
    \item \textbf{Construction of the reduced order system:} This involves matrix products arising when substituting \eqref{modified augmented basis S} into \eqref{PDE} for $\theta=\theta_\alpha$.
\end{itemize}

\subsubsection*{Online Computational Cost}
The online phase includes the computation of the reduced-order system solution and the reconstruction of the numerical solution. Specifically:

\begin{itemize}
    \item \textbf{Computation of the coefficients $\{a_i\}_{i=1}^N$ in the representation \eqref{modified augmented basis S}:} This step requires solving the ROM, which is an ordinary differential system derived from substituting the field representation \eqref{modified augmented basis S} into \eqref{PDE} for the updated parameter values $\theta_\alpha = \theta_0 + \alpha \dtheta$. The computational cost for this task depends on the selected methodology for solving the ODE. 
    \item \textbf{Reconstruction of the solution $u_{\theta_\alpha}$ \eqref{modified augmented basis S}:} This task involves computing scalar vector products to reconstruct the solution.
\end{itemize}

\section{Application to Two-Dimensional Seismic Wave Problems}

In this section, we demonstrate how to apply MOR-T$_L$ to a specific problem, selecting a class of 2D seismic wave problems as our test case. The choice of seismic wave stems from our interest in monitoring CO$_2$ injection in subsurface environments, a process that involves solving an inverse problem where the goal is to reconstruct the velocity profile of the medium. Efficiently solving the corresponding direct problem is a crucial step in achieving this objective.
In many inversion processes, such as Full Waveform Inversion (FWI) (see, e.g., \cite{virieux_overview_2009, metivier_full_2013, hu_feasibility_2022}), it is necessary to solve the same direct problem multiple times, but with different velocity models. These models typically represent perturbations of the original medium, which are determined by solving an adjoint problem. As such, each iteration of the inversion requires solving multiple direct wave problems, often with fine discretization meshes to ensure high accuracy -- especially important when solving ill-posed problems. This results in large-scale, computationally expensive problems, which can render the inversion process computationally infeasible.\\
We propose applying MOR-T$_L$ to this context, demonstrating that a carefully designed ROM can significantly reduce computational costs while maintaining accuracy. This, in turn, can restore efficiency to inversion procedures like FWI.

\subsection{The problem prototype}
Here, we define the boundary value problem for the acoustic wave problem, followed by the set of problems that characterize the Fréchet derivatives of the solutions. 
The initial value problem we consider is given by:
\begin{equation}\label{eq:acoustic wave equation}
    {\bf (IVP)}
    \begin{cases}
        \displaystyle \theta \frac{\partial^2 u_\theta}{\partial t^2}(\x,t)-\Delta{u_\theta(\x,t)} = f(\x,t) &\text{in} \quad \mathbb{R}^2 \times \left(0,T\right) \\ 
        \displaystyle u_\theta(\x,0) = u_0(\x), \; \dfrac{\partial u_\theta}{\partial t}(\x,0) = u_1(\x) &\text{in} \quad \mathbb{R}^2
    \end{cases}
\end{equation}
The parameter $\theta:=\theta(\x)$, $\x \in \mathbb{R}^2$, represents the propagation velocity as it writes $\theta(\x) = \frac{1}{c(\x)^2}$ with $c(\x) > 0$ the functional representing the acoustic wave speed in the domain. The unknown $u_\theta$ is the pressure wave field, $f^{(0)} = f$ is the source term independent of the parameter $\theta$.
Then, $v_\theta^\ell$, the Fréchet derivative of order $\ell$ in the direction of $\dtheta$, satisfies the following collection of initial value problems
\begin{equation}\label{eq:acoustic frechet}
    \begin{cases}
        \displaystyle \theta \frac{\partial^2 v_\theta^\ell}{\partial t^2} - \Delta{v_\theta^\ell} = f^{(\ell)} &\text{in} \quad  \mathbb{R}^2 \times (0,T)\\
        \displaystyle v_\theta^\ell(\x,0) = 0, \quad \frac{\partial v_\theta^\ell}{\partial t}(\x,0) = 0 &\text{in} \quad  \mathbb{R}^2 \\
    \end{cases}
\end{equation}
with $f^{(0)} = f$ and $\displaystyle f^{(\ell)} = - \ell \frac{\partial^2 v_\theta^{\ell-1}}{\partial t^2}$, $1 \leq \ell \leq L$.\\

\begin{remark}
Due to the separable structure of the acoustic wave equation, with the wave operator splitting into space and time components, each right-hand side of the successive Fréchet derivatives depends on the solution from the previous order. As a result, all systems can be solved sequentially at each time step.
\end{remark}

\subsection{The formulation in a Bounded Domain}

To solve the initial value problems \eqref{eq:acoustic wave equation} and \eqref{eq:acoustic frechet} we use the SEM \cite{patera_spectral_1984, komatitsch_spectral_1998}. SEM has the advantage of giving rise to a diagonal mass matrix which in turns can easily accommodate the use of explicit time discretization scheme. This is a clear advantage for large scale wave applications, hence helping in controlling an intensive use of memory. \\
As with any domain-based discretization scheme, the problem must first be reformulated within a bounded domain $\Omega$. This is achived by introducing an artificial external boundary, to which a boundary condition is applied. Therefore, the boundary value problems considered in the bounded domain are given by:

\begin{equation}\label{BVP acoustic vl}
    \begin{cases}
        \displaystyle \theta \frac{\partial^2 v_\theta^\ell}{\partial t^2} - \Delta{v_\theta^\ell} = - \ell \frac{\partial^2 v_\theta^{\ell-1}}{\partial t^2} &\text{in} \quad  \Omega \times (0,T)\\
        \displaystyle B_\theta v_\theta^\ell = 0 &\text{on} \quad \partial \Omega \times (0,T) \\
        \displaystyle v_\theta^\ell(\x,0) = 0, \quad \frac{\partial v_\theta^\ell}{\partial t}(\x,0) = 0 &\text{in} \quad  \Omega \\
        \displaystyle v_\theta^0 = u_\theta 
    \end{cases}
\end{equation}

We have chosen two boundary conditions for this formulation. The first is the Dirichlet boundary condition, where $B_\theta$ is the identity operator, chosen for its simplicity in implementation. The second condition is a first order Absorbing Boundary Condition (ABC) \cite{givoli_high-order_2004, engquist_absorbing_1977} where the operator $B_\theta$ depends on the problem's parameter $\theta$. This choice enables us to assess the efficacy and robustness of the reduced basis with respect to changes in the boundary condition resulting from parameter variations. We aim to verify that MOR-T$_L$ continues to provide a solution with an accuracy level comparable to that obtained by solving the full problem using SEM when both methods are equipped with the same exterior boundary condition. It is important to emphasize that when solving practical or more realistic problems, more accurate boundary conditions would be necessary.

\subsection{The basis functions}
In this section, we describe the practical construction of basis functions used in MOR-T$_L$. It is common in reduced-order methods for the basis functions to be computed using a truncated SVD of the snapshot matrices $S^\ell$, given by equation \eqref{Sl}. This approach is based on POD \cite{bergmann_optimal_2008, stefanescu_poddeim_2013, bui-thanh_model_2008, hawkins_model_2024}. However, there are several limitations to this method, with the following three being particularly noteworthy: (a) the lack of a rigorous criterion - or at least practical guidelines - for selecting the snapshots in $S^\ell$, (b) the potential for $S^\ell$ to become too large to store in memory for large-scale systems, and (c) even when can be stored, the decomposition process itself can become computationally excessive and ultimately prohibitive. 

Hence, we propose a strategy where the SVD is replaced with a truncated QR decomposition based on a Gram-Schmidt (GS) orthonormalization process \cite{leon_gramschmidt_2013} allowing us to compute the basis as we solve \eqref{eq:acoustic wave equation} without storing the full matrices $S^\ell$. This is the purpose of the following.\\

Let $\bfs^\ell_1$ be the snapshot of the (nonzero) solution at a time $t^\ell_1$ given by \eqref{snapshot vector sl}. Let $\mphi^\ell_1$ be the first basis vector defined as:
\begin{equation*}
    \begin{aligned}
        \mphi^\ell_1 &= \frac{\bfs^\ell_1}{\lVert \bfs^\ell_1 \rVert}; \quad \ell=0,\cdots,L  
    \end{aligned}
\end{equation*}
where $\lVert . \rVert$ should be understood as the discrete norm associated with the discrete chosen inner product $\scalp{.}{.}$.\\
We compute the $i_\ell^{th}$ basis functions $i_\ell>1$ in terms of the solution $\bfs^\ell_j$, $j=1,\cdots,N_t$, and the previous basis vectors $\mphi^\ell_k, k=1,\cdots,i_\ell-1$ as:
\begin{equation}
    \begin{aligned}\label{QR algo}
        \Tilde{\mphi}^\ell_{i_\ell} &= \bfs^\ell_j - \sum_{k=1}^{i_\ell-1}\scalp{\bfs^\ell_j}{\mphi^\ell_k} \mphi^\ell_k; \quad \ell=0,\cdots,L.
    \end{aligned}
\end{equation}
If the following condition is satisfied
%The truncation of this basis is done through a criterion on the norm of $\Tilde{\mphi}^\ell_{i_\ell}$ compared to the norm of $\bfs^\ell_j$, follows the rule: for a given $\epsilon>0$, if
\begin{equation}\label{truncation criterion}
    \begin{aligned}
        \lVert\Tilde{\mphi}^\ell_{i_\ell}\rVert^2 > \epsilon \lVert \bfs^\ell_j \rVert^2; \quad \ell=0,\cdots,L 
    \end{aligned}
\end{equation}
for a prescribed value for the positive number $\epsilon$, we proceed to the normalization of the vector.
\begin{equation}\label{QR normalization}
    \mphi^\ell_{i_\ell} = \frac{\Tilde{\mphi}^\ell_{i_\ell}}{\lVert\Tilde{\mphi}^\ell_{i_\ell}\rVert}; \quad \ell=0,\cdots,L. 
\end{equation}
If not, the basis function is not selected. We re-iterate the process (\eqref{QR algo} and \eqref{truncation criterion}) and compute $\Tilde{\mphi}^\ell_{i_\ell}$ with $\bfs^\ell_{j+1}$.\\
The stopping criterion \eqref{truncation criterion} gives us a quantitative way to select the snapshots.
However, GS process is subject to instability error for a large amount of data and may bias this criterion. To reduce this effect, re-orthogonalization must be applied every few selected snapshots. At the end of the process/simulation, we have constructed a basis $\{\mphi^\ell_i\}_{i=1}^{M_\ell}$, with $M_\ell \leq N_t$. In practice, $M_\ell$ depends on the source considered $f$ as well as on the values of the parameter $\theta_0$.\\

To create the final basis, we gather all the modified snapshot matrices $\Tilde{S}^\ell$ in the augmented snapshot matrix $S^{\mathrm{aug}}$ given by \eqref{augmented S}. Then, before applying the reduction order technique on $S^{\mathrm{aug}}$, we proceed to a re-ordering of the column vectors by increasing time, that is, we start from the vector $\mphi^0_1$ obtained at time $t^0_1$, the next vectors are all the vectors $\mphi^\ell_1$ obtained at time $t^\ell_1$, if any, Then we consider the vectors obtained at $t^\ell_2$, if any, starting from $\ell = 0$ up to $\ell = L$, etc. We then obtain after the completion of Gram-Schmidt process a modified full-rank augmented snapshot matrix $\Tilde{S}^{\mathrm{aug}}$ given by \eqref{modified augmented basis S} with our final reduced basis being $\{\mphi_i\}_{i=1}^{N}$, $ N \leq \displaystyle \sum_{\ell=0}^{L} M_\ell$.\\

\begin{remark} The following points are noteworthy:\\
%\begin{itemize}
{\bf{a.}} In addition to being memory-efficient compared to the SVD approach, the QR strategy is easy to parallelize, as each scalar product can be computed independently over subdomains. This parallelization accelerates the process and helps reduce numerical instabilities in the scalar products.\\
{\bf{b.}} The method's efficiency depends on the choice of $\epsilon$. If $\epsilon$ is too small, it may lead to the selection of an excessively large number of basis functions, increasing computation time and potentially causing numerical instabilities. On the other hand, if $\epsilon$ is too large, the number of basis functions may be insufficient, resulting in an overly limited approximation subspace. Therefore, selecting $\epsilon$ requires a balance. In practice, we find that setting $\epsilon=0.01$ offers a good compromise, yielding excellent results.
%\end{itemize}
\end{remark}

\subsection{\textbf{The resulting reduced order model}}

We now introduce the reduced-order model, represented by the boundary value problem

\begin{equation}\label{pod wave}
    \begin{cases}
        \displaystyle \sum_{i=1}^N \theta\frac{\partial^2 a_i(t)}{\partial t^2}\phi_i(\x)
            - a_i(t)\Delta{\phi_i(\x)}  
            = f(\x,t) &\text{in} \quad  \Omega \times ]0,T[ \\
        \displaystyle B_\theta \sum_{i=1}^N a_i(t) \phi_i(\x)
            = 0 &\text{on} \quad  \partial \Omega \times ]0,T[ \\
        \displaystyle a_i(0) = \scalp{u_0(\x)}{\phi_i(\x)}, \; \dfrac{\partial a_i}{\partial t}(0) = \scalp{u_1(\x)}{\phi_i(\x)} &\text{in} \quad \Omega
    \end{cases}
\end{equation}
resulting from substituting \eqref{modified augmented basis S} into \eqref{eq:acoustic wave equation}.\\
Since $u_{\theta_0}(\x,t)$ is supposed to be obtained through a SEM formulation, then the $\{\phi_i(\x)\}_{i=1}^N$ can be written in the SEM basis as $\phi_i(\x) = \displaystyle \sum_{k=1}^{N_x} \phi_{k,i} \psi_k(\x)$. Considering the case where $B_\theta=\sqrt{\theta}\frac{\partial }{\partial t} + \frac{\partial }{\partial \bf{n}}$, and gathering all the $\phi_k,i$ in a matrix $\mphi$, the $f_k$ in a vector $\mathbf{f}$, and the $a_i(t)$ in a vector $\bfa(t)$, $\forall k=1,\cdots,N_x$, $\forall i=1,\cdots,N$, one can write after introducing a weak formulation of \eqref{pod wave}:

\begin{equation}\label{eq:semi discrete pod acoustic}
    \begin{aligned}
        \mphi^T\mathbf{M_\theta}\mphi \frac{\partial^2 \bfa(t)}{\partial t^2}
            +  \mphi^T\mathbf{D_{\sqrt{\theta}}}\mphi \frac{\partial \bfa(t)}{\partial t}
            +  \mphi^T\mathbf{K}\mphi \bfa(t)
            &= \mphi^T \mathbf{f}(t)
    \end{aligned}
\end{equation}
where $\mathbf{M}_\theta$, $\mathbf{D}_{\sqrt{\theta}}$, and $\mathbf{K}$ are the mass, the damping, and the stiffness matrices whose respective entries are: 

\begin{equation}
    \begin{aligned}
        \int_{\Omega} \theta \psi_k(\x)\psi_l(\x) d\x; \quad 
        \int_{\partial\Omega} \sqrt{\theta} \psi_k(\x) \psi_l(\x) d\mathbf{S}; \quad
        \int_{\Omega} \nabla{\psi_k(\x)}.\nabla{\psi_l(\x)} d\x
    \end{aligned}
\end{equation}

For the time discretization we choose a leap-frog scheme for simplicity, but any other suitable scheme may be used. We set $\mathbf{M}_{\mphi} = \mphi^T\mathbf{M_\theta}\mphi$, $\mathbf{D}_{\mphi} = \mphi^T \mathbf{D_{\sqrt{\theta}}} \mphi$, $\mathbf{K}_{\mphi} = \mphi^T\mathbf{K}\mphi$ and $\mathbf{f}_{\mphi} = \mphi^T f$. Therefore, for $n = 1,\cdots,N_t -1$,

\begin{equation}\label{POD formulation abc}
    \begin{aligned}
            \bfa^{n+1} 
            = \left( \mathbf{M}_{\mphi} + \frac{\Delta t}{2} \mathbf{D}_{\mphi}\right)^{-1}
            \left( 
                \left( 2\mathbf{M}_{\mphi} - \Delta t^2 \mathbf{K}_{\mphi} \right) \bfa^n
                - \left( \mathbf{M}_{\mphi} - \frac{\Delta t}{2} \mathbf{D}_{\mphi} \right) \bfa^{n-1}
                + \Delta t^2 \mathbf{f}_{\mphi}
            \right).
    \end{aligned}
\end{equation}

\begin{remark}
In the case of a homogeneous Dirichlet boundary condition, that is, $\displaystyle \sum_{i=1}^N a_i(t) \phi_i(\x) = 0$, $\forall \x \in \partial \Omega$, $\mathbf{D}_{\sqrt{\theta}}=0$. Then, modulo a change in the stiffness matrix and the source term, we obtain:
\begin{equation}\label{POD formulation dirichlet}
    \begin{aligned}
            \bfa^{n+1} 
            = 2 \bfa^n - \bfa^{n-1} - \Delta t^2 \Mphi^{-1} \left( \Kphi \bfa^n - \mathbf{f}_{\mphi} \right); \quad n = 1,\cdots,N_t -1   
    \end{aligned}
\end{equation}
\end{remark}

\begin{remark}
A judicious choice of inner product during the QR decomposition, as described in Section 4.3, can enhance the computational efficiency of the MOR-T$_L$ strategy. For the application under consideration, we select the inner product associated with the stiffness matrix $\mathbf{K}$ , defined as $\scalp{u}{v}_\mathbf{K} = \bfu^T \mathbf{K} \bfv$. This choice offers two key computational advantages: (a) due to the leap-frog time-stepping scheme, only half of the inner product is inherently computed at each time step when solving the resulting SEM linear system, and (b) it eliminates the need to compute the projections $\mphi^T \mathbf{K} \mphi$, as the basis is orthonormal with respect to $\mathbf{K}$ by construction, thereby resulting in significant computational cost savings. The only projections that need to be computed, as indicated in (4.7), are those associated with the mass $\mathbf{M}_\theta$ and damping $\mathbf{D}_{\sqrt{\theta}}$ matrices. These projections are much less computationally demanding, since both the mass and damping matrices are diagonal, a property arising from the SEM formulation.
\end{remark}

\section{Assessing cost-effectiveness and performances\protect\footnote{All these numerical experiments have been performed on a system equipped with an 11th Gen Intel Core i5-1135G7 (4 cores, 8 threads, 2.40GHz), Mesa Intel Xe Graphics (TGL GT2). Python 3.10.12 with Numpy 1.21.5 multi-threaded for BLAS operations}of MOR-T$_L$ v.s. SEM in wave problems}

In this section, we present a numerical investigation aimed at evaluating the cost-effectiveness of the proposed MOR strategy compared to the standard SEM approach for solving the class of two-dimensional wave problems \eqref{eq:acoustic wave equation}. Our analysis is divided into two parts. First, we examine the cost-effectiveness of using QR versus SVD basis representations in the context of MOR. Next, we compare the computational cost of MOR-T$_L$ with that of the SEM method, assessing performance for a range of velocity values and a prescribed accuracy level. We conduct these experiments for both homogeneous and heterogeneous media to capture a broad set of scenarios. Additionally, we investigate the impact of the exterior artificial boundary choice on computational cost, presenting results with both the simplest homogeneous Dirichlet boundary condition and a more realistic absorbing boundary condition.\\
We consider $u(\x,0) = 0$, $\frac{\partial u}{\partial t}(\x,0) = 0$, but a Ricker source of the form:

\begin{equation}\label{ricker 2}
    \begin{aligned}
        f(\x,t) = \left(1-\frac{1}{2}\omega^2 (t-t_0)^2\right)e^{-\frac{1}{4}\omega^2(t-t_0)^2}\delta(\x - \x_s)
    \end{aligned}
\end{equation}
where $\omega = 2\pi f_s$ the pulsation with $f_s$ the central source frequency and $t_0$ a delay.\\

\subsection{Assessing the Cost-Effectiveness of QR and SVD Basis Representations}
\begin{figure}[H]
    \centering
    \includegraphics[width=0.5\textwidth]{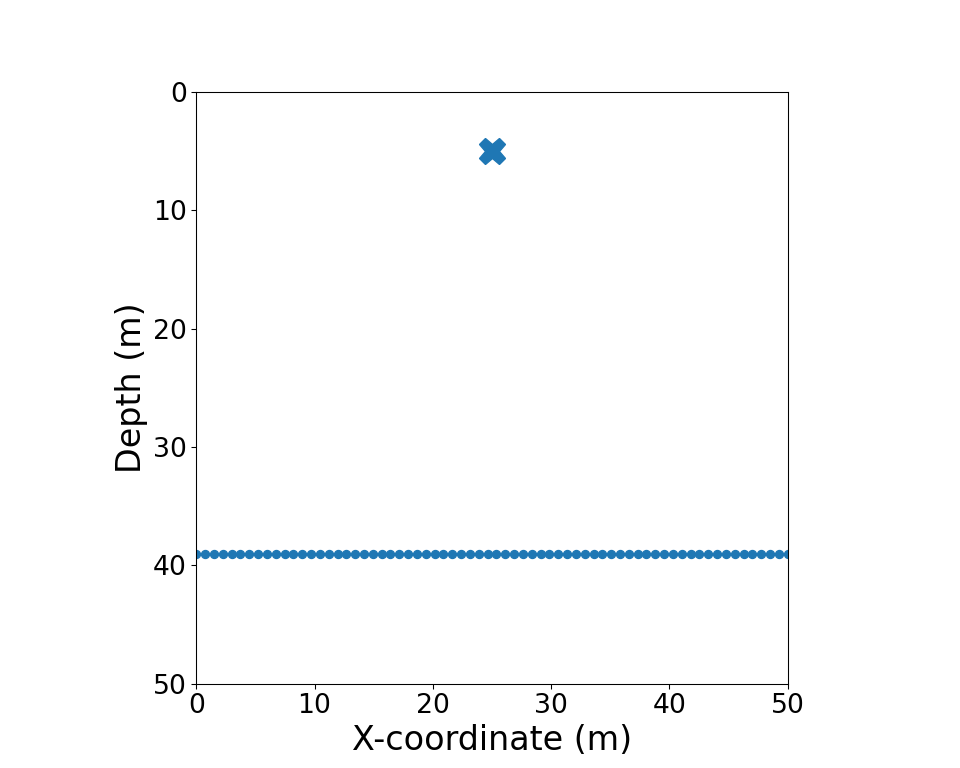} 
    \caption{{\bf{\footnotesize 2D computational domain $\Omega$ showing its depth and the horizontal axis.}} \\ \footnotesize \it The figure illustrates the location of the source f and the positions of the 202 receivers, which are aligned along a horizontal line to measure the pressure field.}
    \label{fig:source receivers}
\end{figure}

In this subsection, we present numerical results comparing the performance and cost-effectiveness of the QR and SVD basis representations. We consider the initial boundary value problem \eqref{eq:acoustic wave equation}, defined over a $50$m $\times 50$m square computational domain $\Omega$, with a source function $f$ represented by the Ricker wavelet \eqref{ricker 2} with a central frequency of $f_s=2.5$Hz and homogeneous initial conditions, i.e., $u_0=u_1=0$. The source is located at $\x_s=(25,5)$ as depicted in Figure \ref{fig:source receivers}. \\%and we have a line of 202 receivers from $\x_{r_1}=(0,39.05)$ to $\x_{r_{202}}=(50,39.05)$ as shown in Figure \ref{fig:source receivers}

The reference solutions are computed using the SEM method with a $201\times201$ Q1 finite element mesh \cite{komatitsch_spectral_1998}. The experiments are run over 4000 time steps with $\Delta t=0.001$, and the model parameter is $\theta=\frac{1}{c^2}$, where $c$ is the velocity in $\Omega$. These experiments are conducted on a fixed, unperturbed velocity, with both homogeneous and heterogeneous media considered. For clarity, results are presented for $\epsilon=0.01$, though results for other values of $\epsilon$ are provided in \cite{jb}.

\subsubsection{Case of Homogeneous Medium}
We begin with a homogeneous medium where the velocity is constant, $c(x)=15$m/s in $\Omega$. Here, we assess the cost-effectiveness of QR and SVD basis representations in two boundary condition scenarios: a simple homogeneous boundary condition and a more realistic absorbing boundary condition.

\subsubsection*{\underline{Case of Homogeneous Dirichlet boundary condition.}}
The results for the homogeneous Dirichlet boundary condition are shown in Figures \ref{fig: svd qr basis functions homogeneous Dirichlet} - \ref{fig:error dirichlet homogeneous QR-SVD} and Table \ref{tab:comparison SVD QR homogeneous dirichlet}. Notable observations include:

\begin{itemize}
    \item Figure \ref{fig: svd qr basis functions homogeneous Dirichlet}, illustrates sample basis functions for both the SVD and QR decompositions. With $\epsilon=0.01$, the QR method yields $N=81$ basis elements. The figure highlights that the QR basis functions capture the finite speed propagation of the wave support, while the SVD functions preserve the frequency content. The QR approach maintains time causality, as each basis function is computed from linear combinations of the solution $u_{\theta_0}(\x,t)$ at the current and previous time steps (see \eqref{QR algo}). 
    \item Figure \ref{fig:traces dirichlet homogeneous SEM QR-SVD} and \ref{fig:error dirichlet homogeneous QR-SVD} show comparisons of the accuracy of the ROM solutions using QR and SVD bases. Time traces (Figure \ref{fig:traces dirichlet homogeneous SEM QR-SVD}) at a selected location $x_r =(39.05,15.42)$ demonstrate excellent agreement with the reference SEM solution for both QR and SVD methods. Similar results were obtained at other source and receiver locations, as detailed in \cite{jb}. Figure \ref{fig:error dirichlet homogeneous QR-SVD} plots the relative error in $L^2$ norm over the computational domain $\Omega$ at every $100^{th}$ time step, confirming that both representations achieve comparable accuracy (below $10^{-2}\%$). Initially, QR significantly outperforms SVD, likely due to the QR method’s first basis function being the exact normalized solution $u_{\theta_0}(\x,t)$ at $t=t_1$. Additionally, results in \cite{jb} confirm that reducing $\epsilon$ improves accuracy, as expected.
    \item Table \ref{tab:comparison SVD QR homogeneous dirichlet} demonstrates the superior computational efficiency of the MOR-T$_0$ formulation. The MOR-T$_0$ method requires approximately 2.5 less computational time to achieve comparable -- or even better -- accuracy than the POD approach.
    \item Final remark: when the source function ceases emission, increasing $\epsilon$ causes the QR algorithm to separate consecutive basis functions by an amount corresponding to a wavelength. This behavior allows significant computational cost savings, as fewer snapshots need to be evaluated, reducing the overall computational load.
\end{itemize}

\begin{minipage}[t]{0.49\linewidth}
\hfill
\begin{minipage}[t]{0.49\linewidth}
    \centering
    \includegraphics[width=\linewidth]{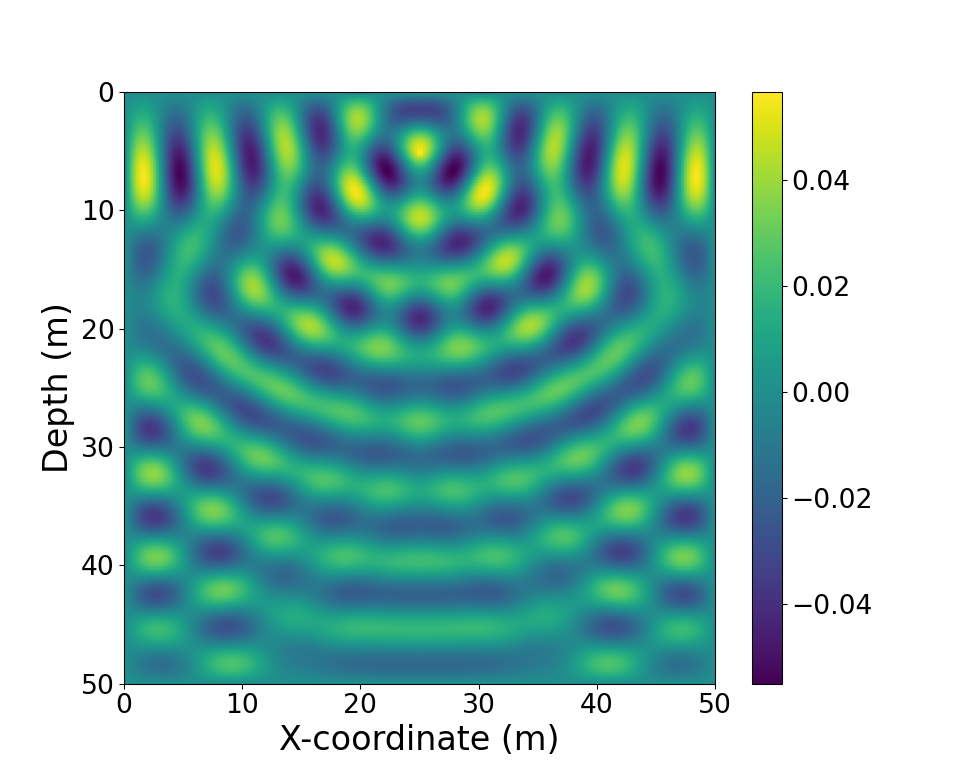}
    \centering{$\phi_1$}
\end{minipage}%
\hfill
\begin{minipage}[t]{0.49\linewidth}
    \centering
    \includegraphics[width=\linewidth]{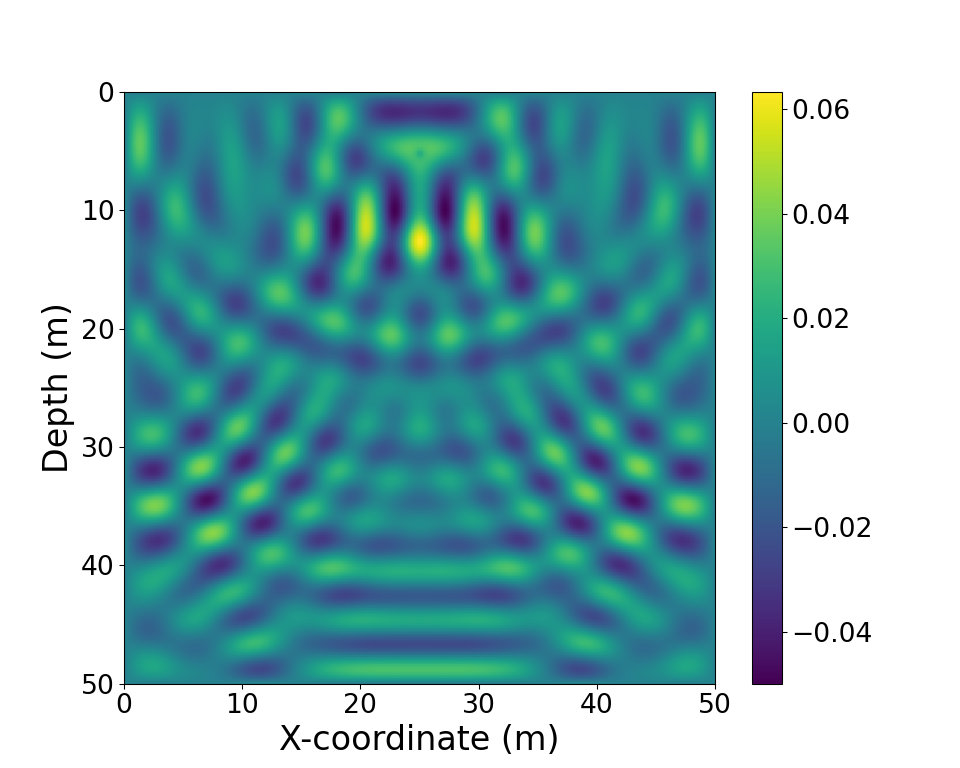}
    \centering{$\phi_{20}$}
\end{minipage}%
\vspace{1em}
\begin{minipage}[t]{0.49\textwidth}
    \centering
    \includegraphics[width=\textwidth]{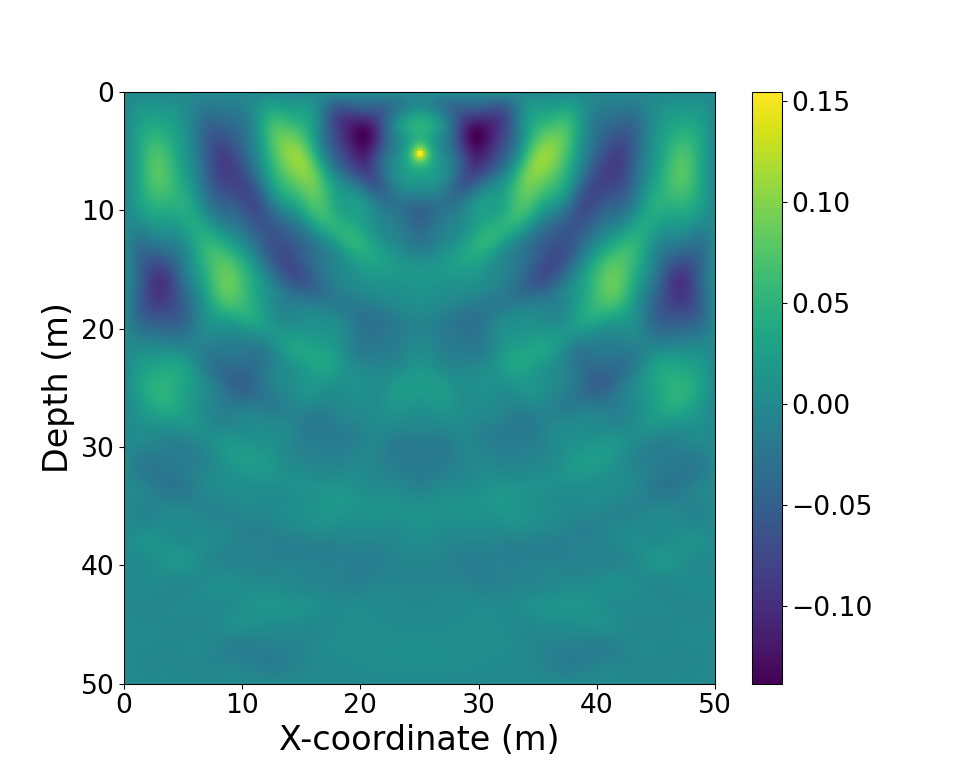}
    \centering{$\phi_{25}$}
\end{minipage}%
\hfill
\begin{minipage}[t]{0.49\textwidth}
    \centering
    \includegraphics[width=\textwidth]{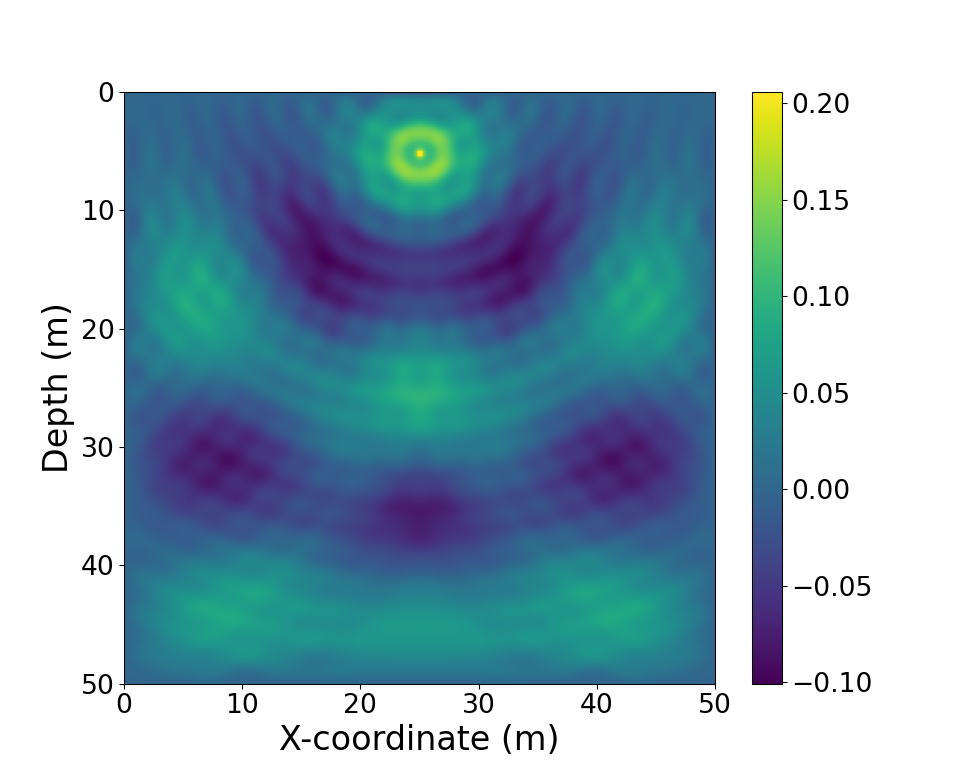}
    \centering{$\phi_{40}$}
\end{minipage}%
\vspace{1em}
\begin{minipage}[t]{0.49\textwidth}
    \centering
    \includegraphics[width=\textwidth]{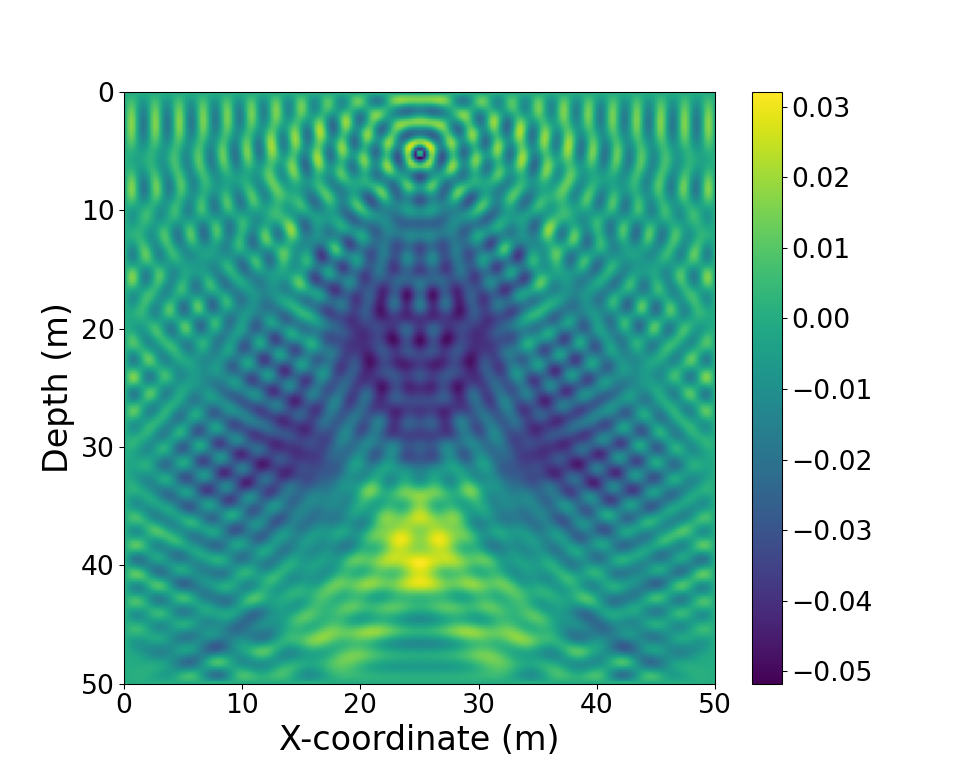}
    \centering{$\phi_{60}$}
\end{minipage}%
\hfill
\begin{minipage}[t]{0.49\textwidth}
    \centering
    \includegraphics[width=\textwidth]{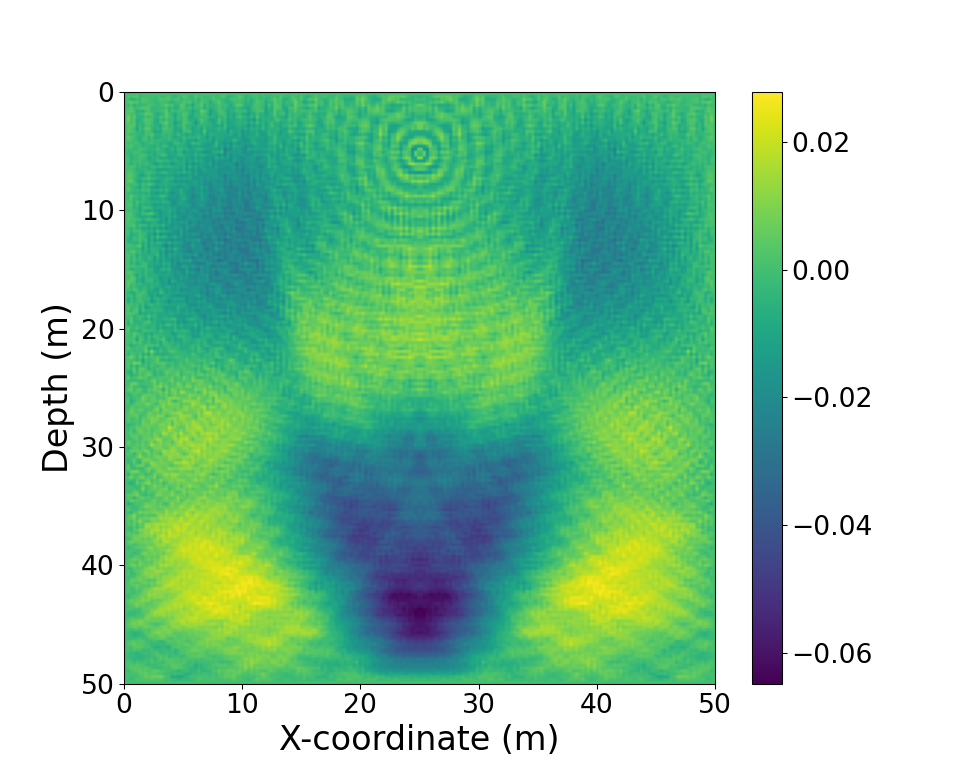}
    \centering{$\phi_{81}$}
\end{minipage}
\centering{SVD}
\end{minipage}
\vrule
\begin{minipage}[t]{0.49\textwidth}
\hfill
\begin{minipage}[t]{0.49\textwidth}
    \centering
    \includegraphics[width=\textwidth]{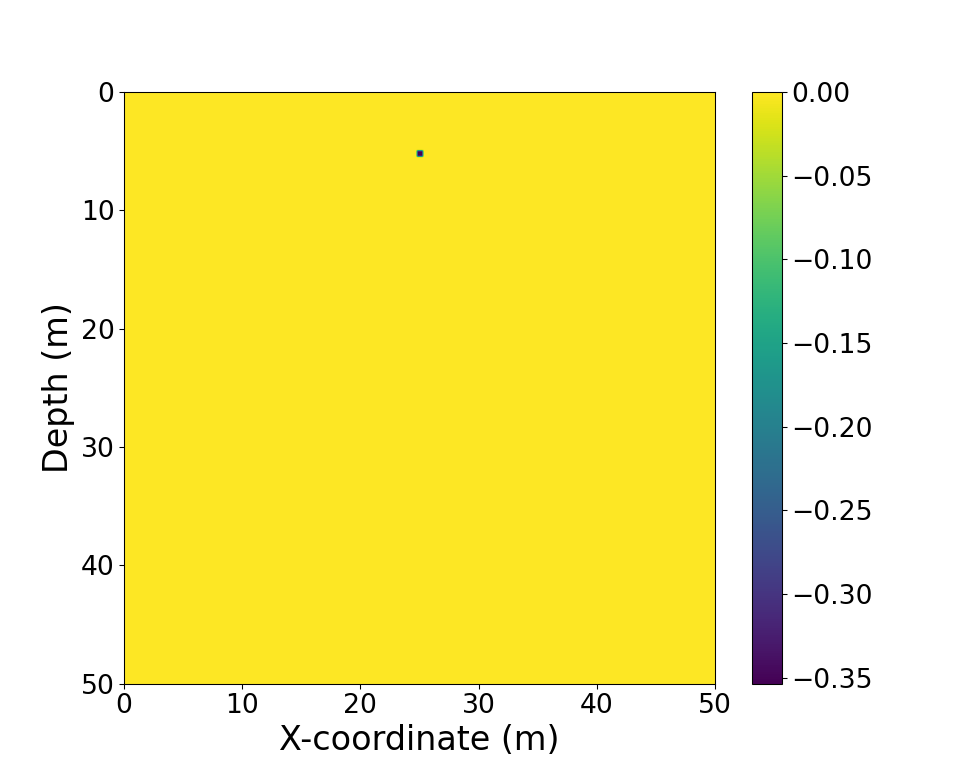}
    \centering{$\phi_1$}
\end{minipage}%
\hfill
\begin{minipage}[t]{0.49\textwidth}
    \centering
    \includegraphics[width=\textwidth]{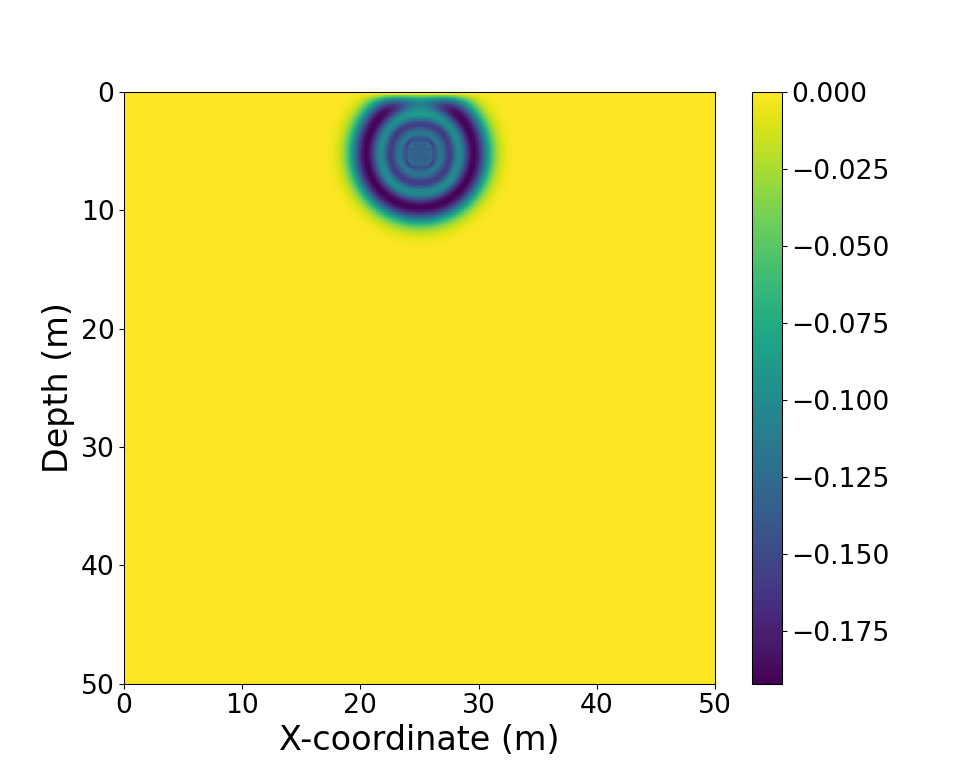}
    \centering{$\phi_{20}$}
\end{minipage}%
\vspace{1em}
\begin{minipage}[t]{0.49\textwidth}
    \centering
    \includegraphics[width=\textwidth]{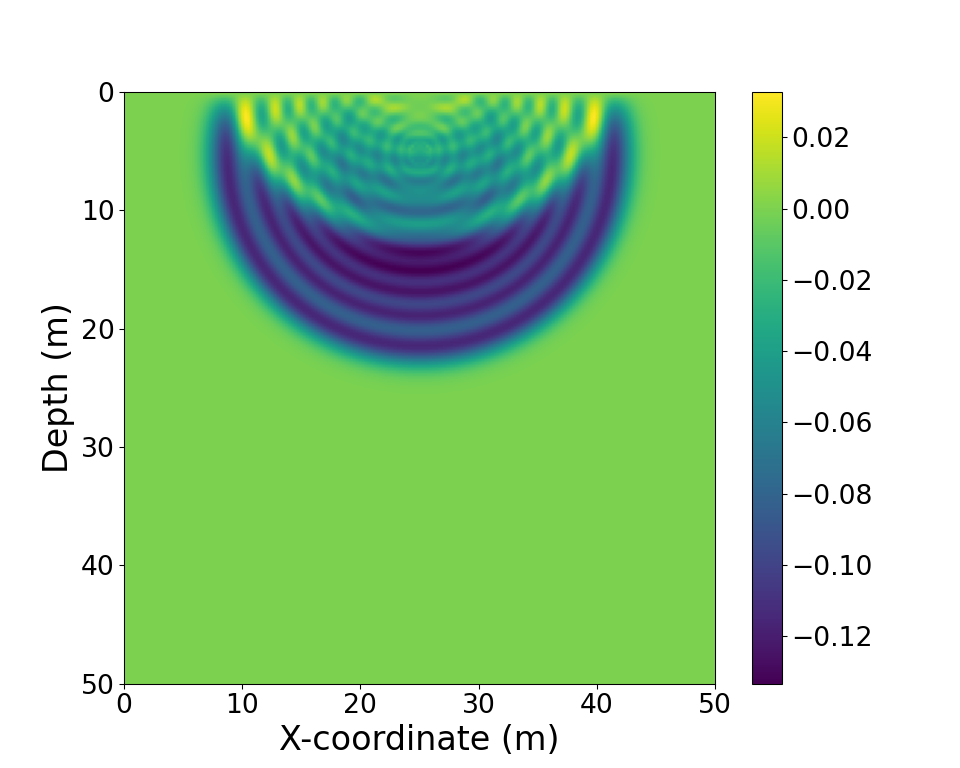}
    \centering{$\phi_{25}$}
\end{minipage}%
\hfill
\begin{minipage}[t]{0.49\textwidth}
    \centering
    \includegraphics[width=\textwidth]{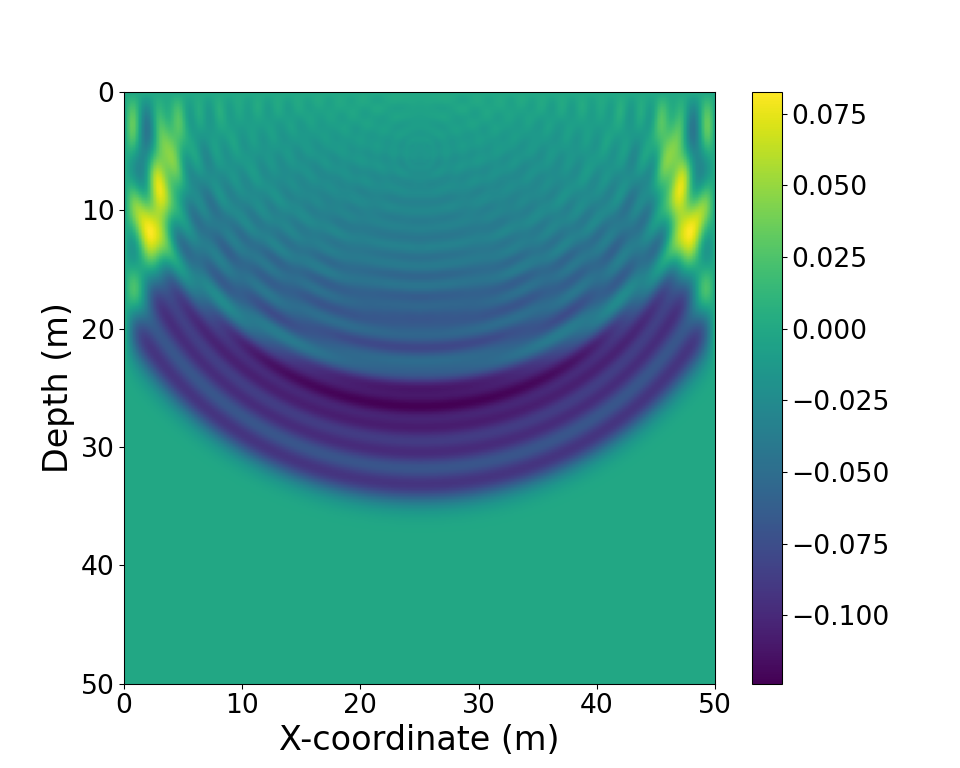}
    \centering{$\phi_{40}$}
\end{minipage}%
\vspace{1em}
\begin{minipage}[t]{0.49\textwidth}
    \centering
    \includegraphics[width=\textwidth]{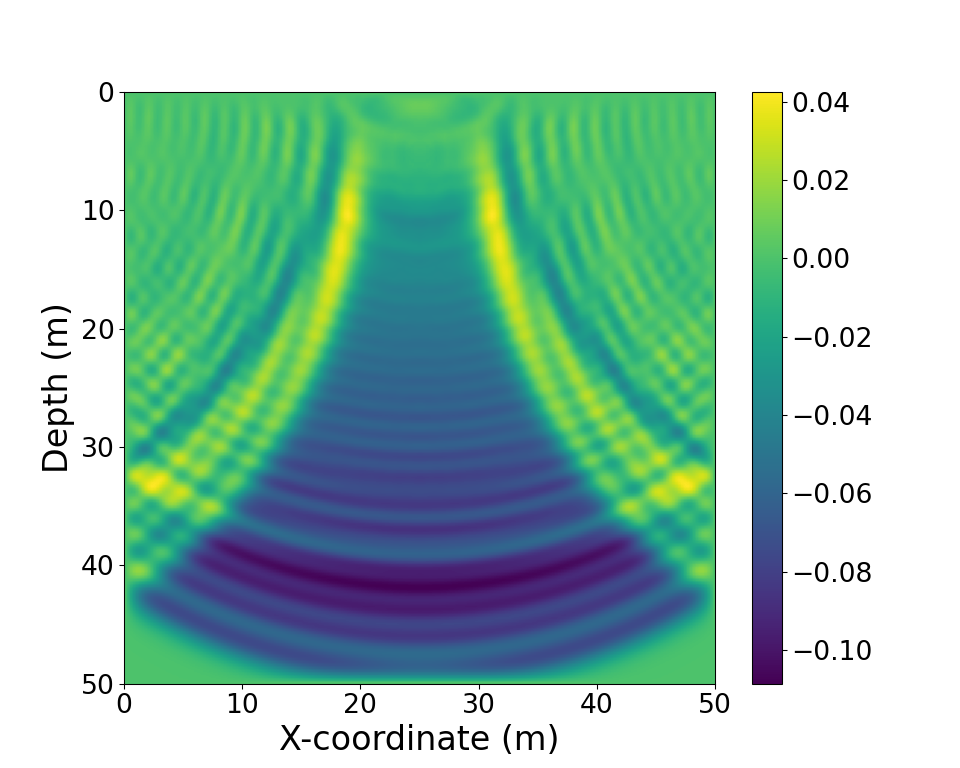}
    \centering{$\phi_{60}$}
\end{minipage}%
\hfill
\begin{minipage}[t]{0.49\textwidth}
    \centering
    \includegraphics[width=\textwidth]{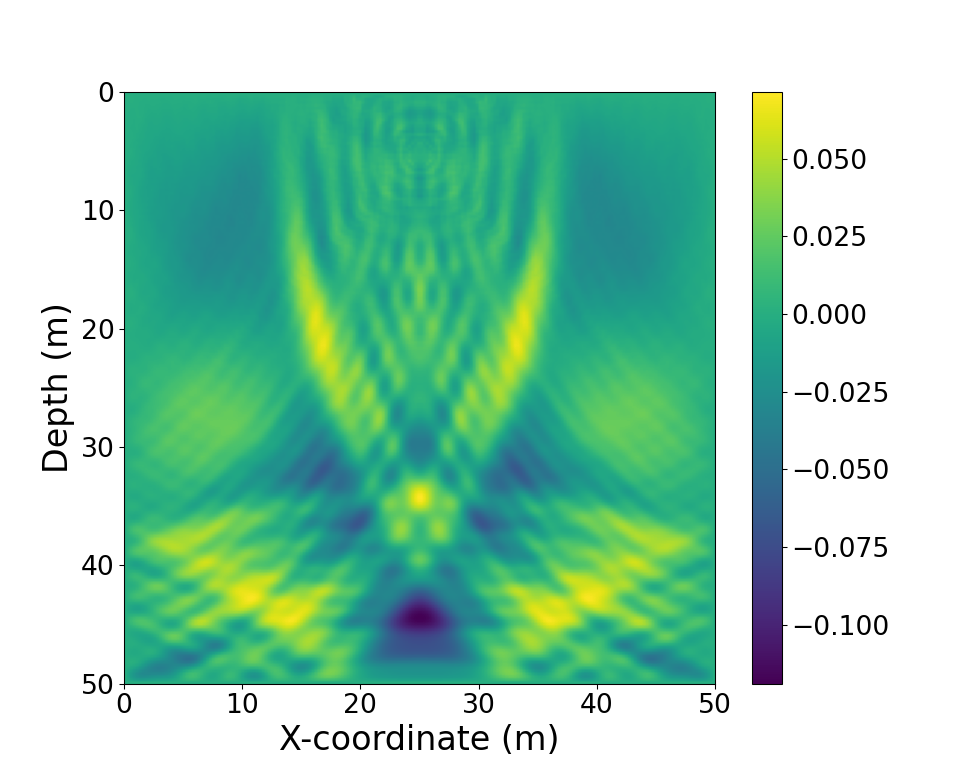}
    \centering{$\phi_{81}$}
\end{minipage}
\centering{QR}
\end{minipage}\\
\captionof{figure}{Sample of basis functions: SVD (left) v.s. QR (right).}
\label{fig: svd qr basis functions homogeneous Dirichlet}

\bigskip 

\begin{minipage}[t]{0.49\textwidth}
\begin{figure}[H]
    \centering
    \includegraphics[width=\textwidth]{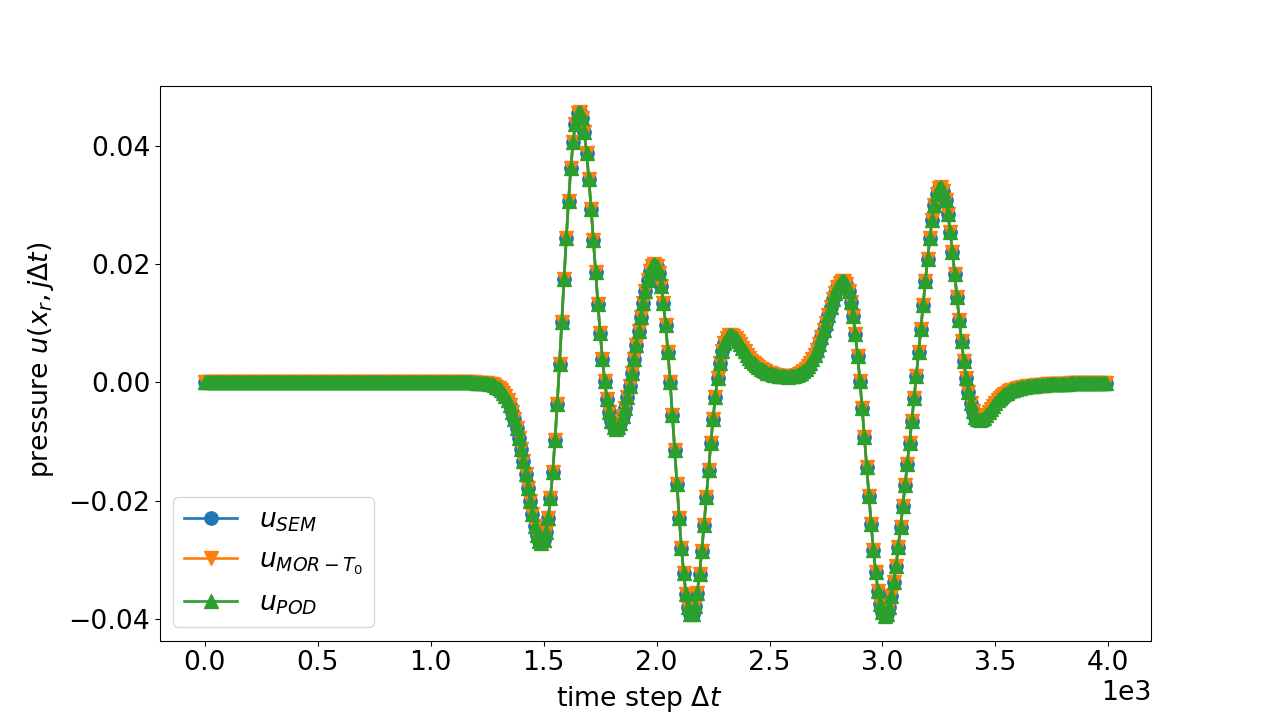} 
    \caption{{\bf{\footnotesize Pressure Field Comparison at $\boldsymbol{x_r = (39.05, 15.42)}$ for a Homogeneous Medium.}} \\ \footnotesize \it Pressure values over $[0,4×10^3]$ are compared between SEM (reference), MOR-T$_0$, and POD for a homogeneous medium with boundary condition $B_\theta= I$.}
    \label{fig:traces dirichlet homogeneous SEM QR-SVD}
\end{figure}
\end{minipage}
\hfill
\begin{minipage}[t]{0.49\textwidth}
\begin{figure}[H]
    \centering
    \includegraphics[width=\textwidth]{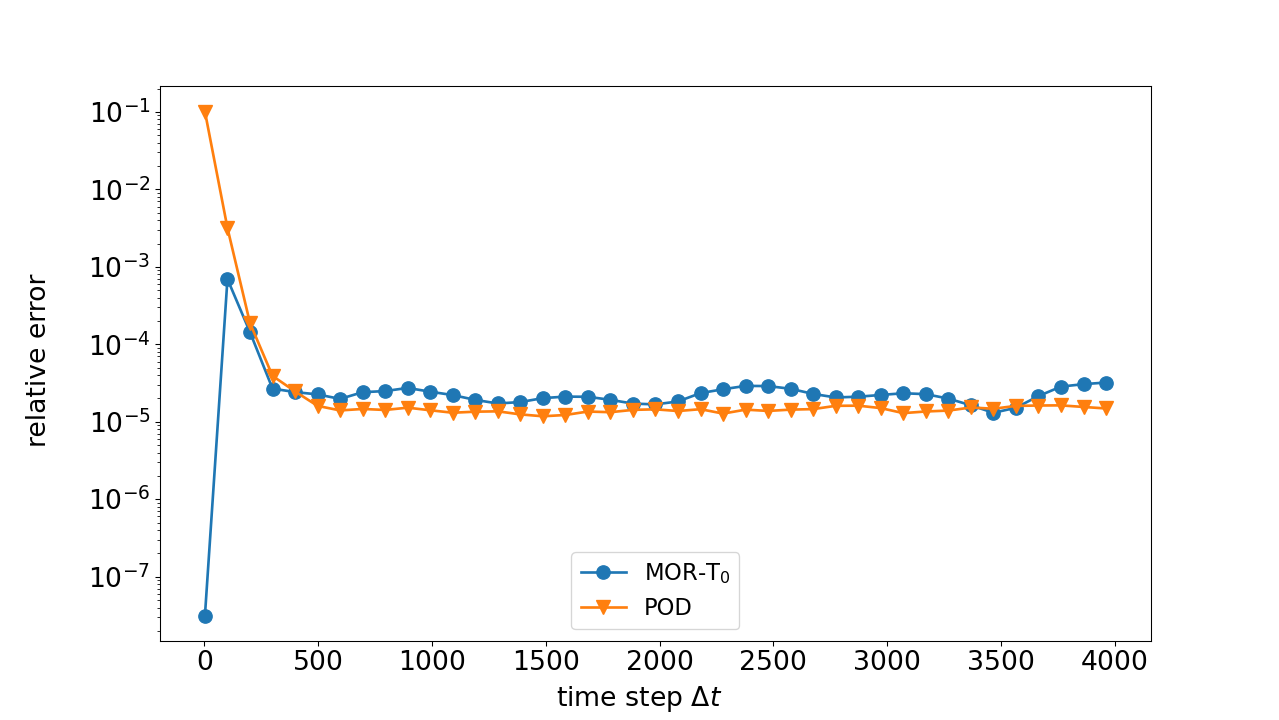} 
    \caption{{\bf{\footnotesize Relative Errors Between MOR-T$_0$, POD, and SEM.}} \\ \footnotesize\it Relative errors over the time interval $[0, 4\times10^3]$ are shown for MOR-T$_0$, and POD compared to SEM. The error is computed as the standard relative L$^2$-norm error across the computational domain $\Omega$ at each time step, for a homogeneous medium $\Omega$ with boundary condition $B_\theta = I$.}
    \label{fig:error dirichlet homogeneous QR-SVD}
\end{figure}
\end{minipage}
\bigskip

\begin{table}[H]
    \centering
    \begin{tabular}{|c|c|c|c|}
        \hline
        Method & Total Computational Time & $\displaystyle \frac{1}{N_t}\sum_{i=1}^{N_t} \frac{\lVert u_{SEM}(\x,t_i) - u_{ROM}(\x,t_i) \rVert_2}{\lVert u_{SEM}(\x,t_i) \rVert_2}$  \\
        \hline
        $SVD$ &150.0s &$6.54 \times 10^{-4}$ \\
        \hline
        $QR$ &66.1s &$5.19 \times 10^{-5}$ \\
        \hline
    \end{tabular}
    \caption{{\bf{\footnotesize Comparison of MOR-T$_0$ and POD: Total Computational Time and Relative Errors.}} \\\footnotesize \it This table presents the total computational time (offline time for constructing the basis and online time for solving the system) and the average relative L$^2$-norm error for MOR-T$_0$ and POD. The errors are computed across the computational domain $\Omega$ and averaged over time steps. These results are obtained in the case of a homogeneous medium $\Omega$ and a homogeneous boundary condition $B_\theta= I$.}
    \label{tab:comparison SVD QR homogeneous dirichlet}
\end{table}

In conclusion, this experiment clearly establishes the MOR-T$_0$ formulation’s superiority over POD in terms of computational efficiency. The QR-based approach also offers strong parallelizability, which can further reduce computation time on parallel platforms. Furthermore, QR decomposition requires less storage since only one snapshot of the solution needs to be stored at any given time, while the SVD approach requires storing the entire matrix $S$. Thus, MOR-T$_0$ outperforms POD both in terms of computational cost and storage requirements -- a crucial consideration for high-frequency wave simulations.

\subsubsection*{\underline{Case of Absorbing Boundary Condition.}}

We now turn to the case where the boundary condition is replaced by an absorbing boundary condition, characterized by the operator $B_\theta=\sqrt{\theta}\frac{\partial }{\partial t} + \frac{\partial }{\partial \bf{n}}$. The results, shown in Figures  \ref{fig:traces abc homogeneous SEM QR-SVD} and \ref{fig:error abc homogeneous QR-SVD}, as well as Table \ref{tab:comparison SVD QR abc}, reveal similar behavior to the Dirichlet case:

\begin{itemize}
    \item Figure \ref{fig:traces abc homogeneous SEM QR-SVD} shows a perfect match between the time trace solutions from both MOR-T$_0$ and POD methods and the reference SEM solution.
    \item Figure \ref{fig:error abc homogeneous QR-SVD}, demonstrates that both methods achieve a comparable accuracy (within $10^{-3}\%$), although the POD method exhibits a noticeable error increase at the initial and final time steps.
    \item Table \ref{tab:comparison SVD QR abc} again confirms the superior computational efficiency of MOR-T$_0$. Specifically, MOR-T$_0$ reduces computational time by approximately 4 times, while providing an accuracy three orders of magnitude better than POD.  
\end{itemize}

In summary, the use of a more sophisticated absorbing boundary condition does not diminish the performance of MOR-T$_0$, maintaining its computational superiority over the POD method.

\begin{minipage}[t]{0.49\textwidth}
\begin{figure}[H]
    \centering
    \includegraphics[width=\textwidth]{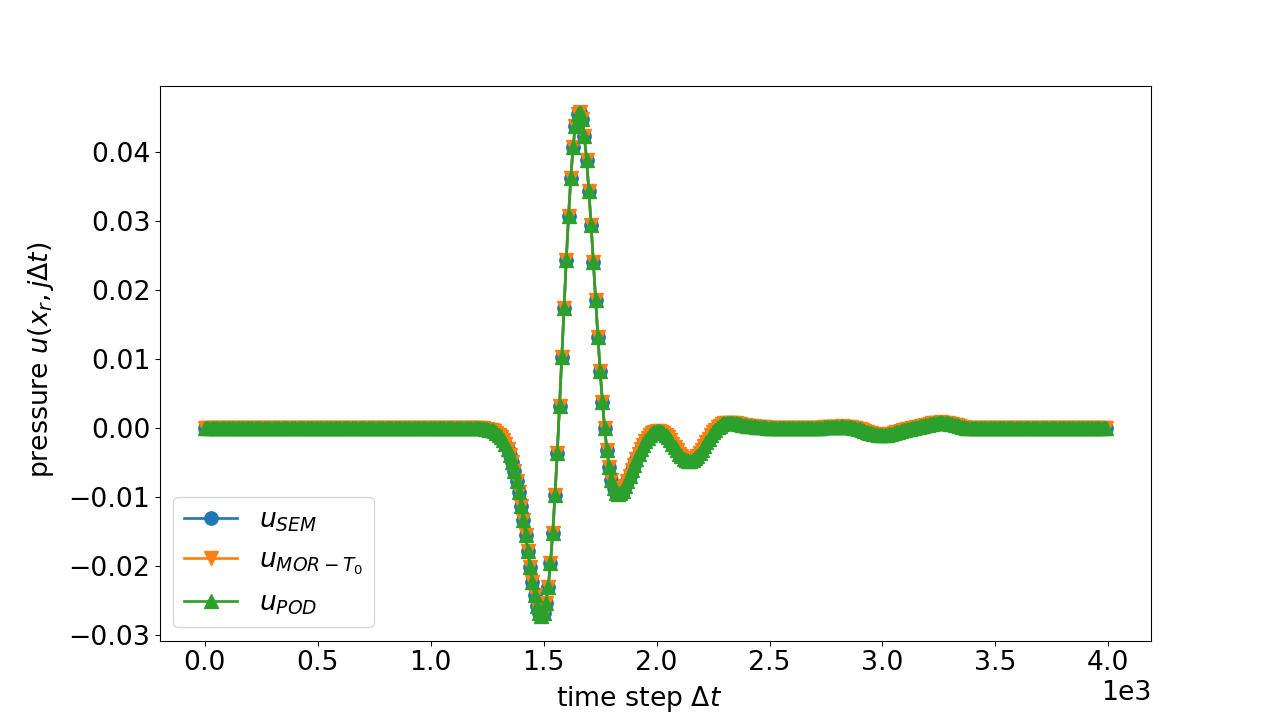} 
    \caption{{\bf{\footnotesize Pressure Field Comparison at \bm{$\x_r = (39.05, 15.42)$} for a Homogeneous Medium.}} \\ \footnotesize \it Pressure values over $[0,4×10^3]$ are compared between SEM (reference), MOR-T$_0$, and POD for a homogeneous medium with boundary condition $B_\theta= \sqrt{\theta}\frac{\partial }{\partial t} + \frac{\partial }{\partial \bf{n}}$.}
    \label{fig:traces abc homogeneous SEM QR-SVD}
\end{figure}
\end{minipage}
\hfill
\begin{minipage}[t]{0.49\textwidth}
\begin{figure}[H]
    \centering
    \includegraphics[width=\textwidth]{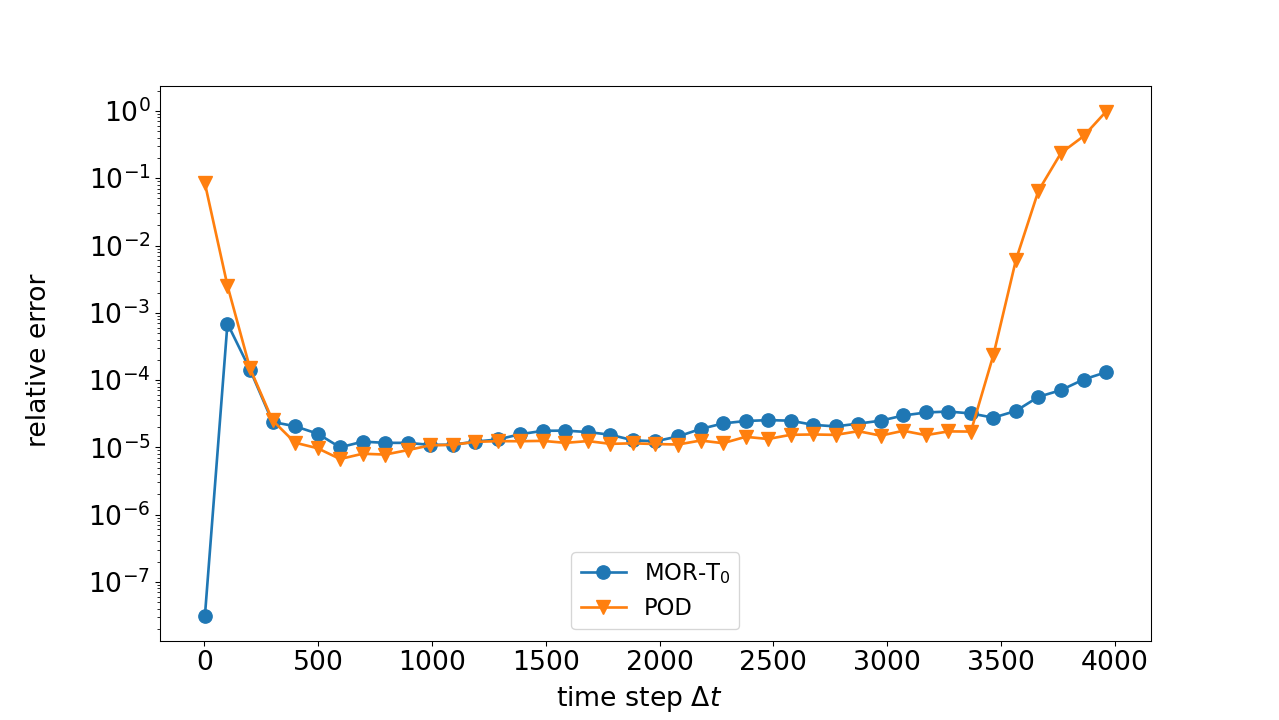} 
    \caption{{\bf{\footnotesize Relative Errors Between MOR-T$_0$, POD, and SEM.}}\\ \footnotesize \it Relative errors over the time interval $[0, 4\times10^3]$ are shown for MOR-T$_0$ and POD compared to SEM. The error is computed as the standard relative L$^2$-norm error across the computational domain $\Omega$ at each time step, for a homogeneous medium $\Omega$ with boundary condition $B_\theta = \sqrt{\theta}\frac{\partial }{\partial t} + \frac{\partial }{\partial \bf{n}}$.}
    \vspace{-0.67cm} 
    \label{fig:error abc homogeneous QR-SVD}
\end{figure}
\end{minipage}

\bigskip
\bigskip

\begin{table}[H]
    \centering
    \begin{tabular}{|c|c|c|c|}
        \hline
        Method & Total Computational Time & $\displaystyle \frac{1}{N_t}\sum_{i=1}^{N_t} \frac{\lVert u_{SEM}(\x,t_i) - u_{ROM}(\x,t_i) \rVert_2}{\lVert u_{SEM}(\x,t_i) \rVert_2}$  \\
        \hline
        $SVD$ &132.2s &$4.04 \times 10^{-2}$ \\
        \hline
        $QR$ &32.4s &$5.56 \times 10^{-5}$ \\
        \hline
    \end{tabular}
    \caption{{\bf{\footnotesize Comparison of MOR-T$_0$ and POD: Total Computational Time and Relative Errors.}}\\ \footnotesize \it This table presents total computational time (offline time for constructing the basis and online time for solving the system) and the average relative L$^2$-norm error for MOR-T$_0$ and POD. The errors are computed across the computational domain $\Omega$ and averaged over time steps. These results are obtained in the case of a homogeneous medium $\Omega$ and a homogeneous boundary condition $B_\theta = \sqrt{\theta}\frac{\partial }{\partial t} + \frac{\partial }{\partial \bf{n}}$.}
    \label{tab:comparison SVD QR abc}
\end{table}

\subsubsection{Case of Heterogeneous Medium}

\begin{figure}[H]
    \centering
    \includegraphics[width=.5\textwidth]{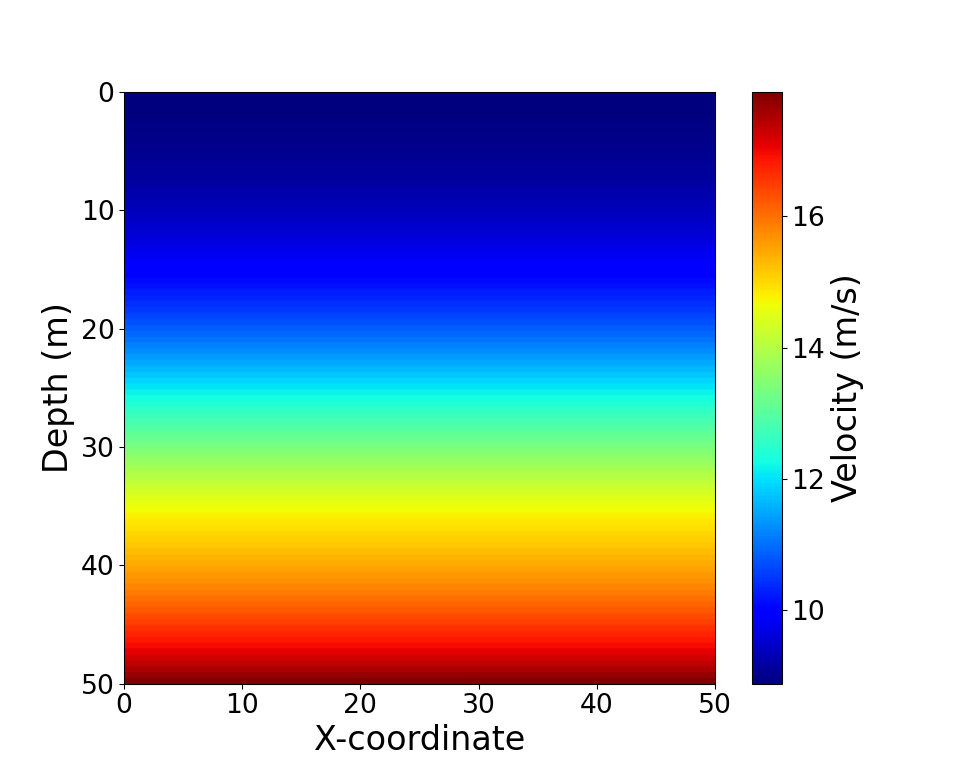} 
    \caption{{\bf{\footnotesize Velocity Distribution in a Heterogeneous Medium.}} \\\footnotesize \it The figure shows the velocity $c_0 = \frac{1}{\sqrt{\theta_0}}$ in the heterogeneous medium $\Omega$, with depth represented on the y-axis and the x-coordinate on the x-axis. The color map indicates the velocity values.}
    \label{fig:model m0 marmousi}
\end{figure}

For the heterogeneous medium, we consider a smooth velocity field $c(\x)$, as shown in Figure \ref{fig:model m0 marmousi}. Results are presented for the absorbing boundary condition, with further details available in \cite{jb}. The time trace solutions and accuracy levels follow the same pattern as the homogeneous case, with perfect agreement and impressive accuracy (below $10^{-2} \%$), as depicted in Figures \ref{fig:traces abc heterogeneous SEM QR-SVD} and \ref{fig:error abc heterogeneous QR-SVD}. This further supports the robustness of the MOR-T$_0$ method across different medium configurations. %Based on Remark \ref{sampling qr}, we selected a sampling strategy where snapshots are taken at a time step equal to one-fifth of the smallest wavelength. Consequently, our algorithm is applied only to snapshots separated by this time step, significantly reducing the computational cost of constructing the reduced basis.

\begin{minipage}{0.49\textwidth}
\begin{figure}[H]
    \centering
    \includegraphics[width=\textwidth]{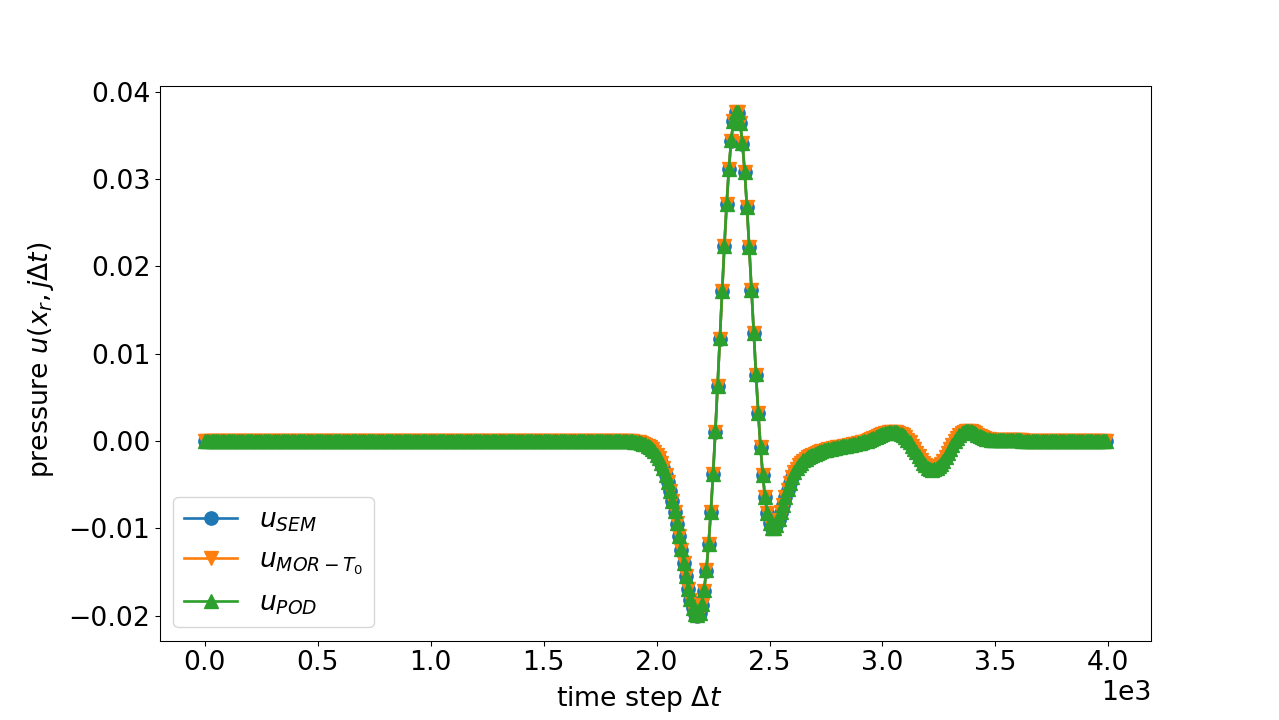} 
    \caption{{\bf{\footnotesize Pressure Field Comparison at \bm{$\x_r = (39.05, 15.42)$} for a Heterogeneous Medium.}} \\\footnotesize \it Pressure values over $[0,4×10^3]$ are compared between SEM (reference), MOR-T$_0$, and POD for a heterogeneous medium with boundary condition $B_\theta= \sqrt{\theta}\frac{\partial }{\partial t} + \frac{\partial }{\partial \bf{n}}$}
    \label{fig:traces abc heterogeneous SEM QR-SVD}
\end{figure}
\end{minipage}
\hfill
\begin{minipage}{0.49\textwidth}
\begin{figure}[H]
    \centering
    \includegraphics[width=\textwidth]{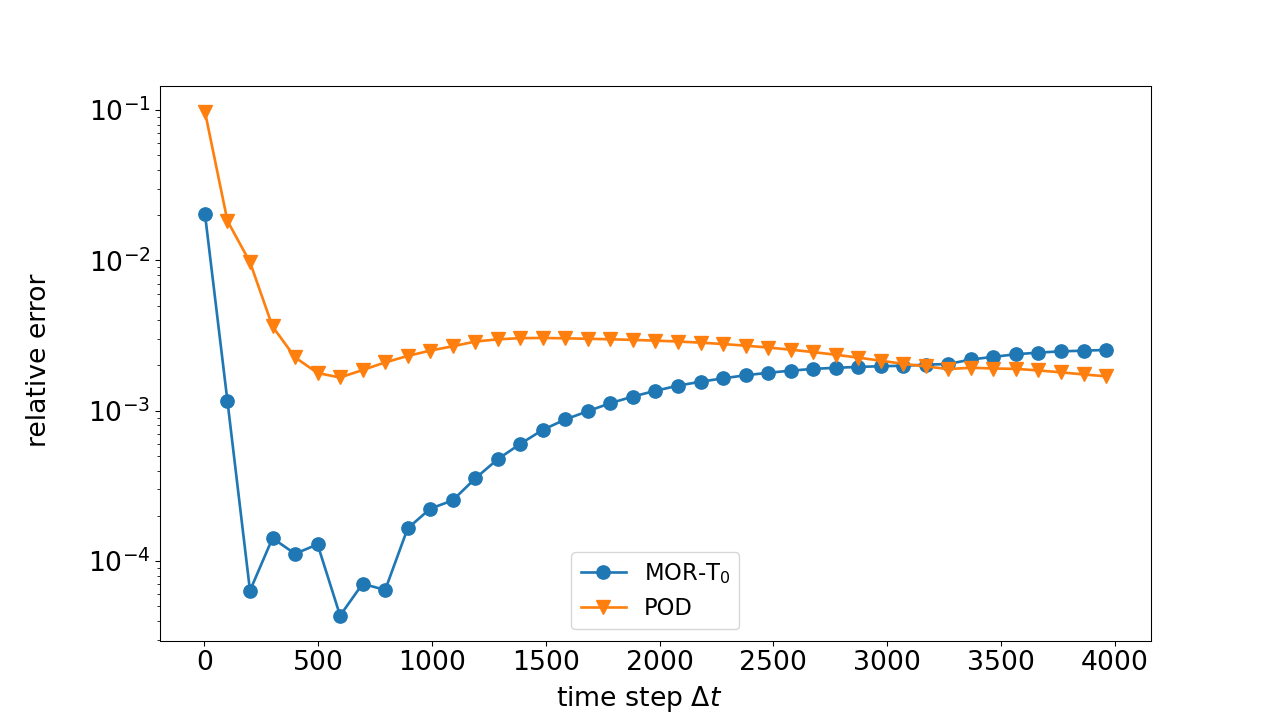}
    \caption{{\bf{\footnotesize Relative Errors Between MOR-T$_0$, POD, and SEM.}}\\ \footnotesize \it Relative errors over the time interval $[0, 4\times10^3]$ are shown for MOR-T$_0$ and POD compared to SEM. The error is computed as the standard relative L$^2$-norm error across the computational domain $\Omega$ at each time step, for a heterogeneous medium $\Omega$ with boundary condition $B_\theta = \sqrt{\theta}\frac{\partial }{\partial t} + \frac{\partial }{\partial \bf{n}}$.}
    \vspace{-0.67cm}
    \label{fig:error abc heterogeneous QR-SVD}
\end{figure}
\end{minipage}
\bigskip
\bigskip

\subsection{Comparative Analysis of MOR-T$_L$ and SEM in Line Search Minimization}
In this subsection, we evaluate both the accuracy and computational efficiency of the proposed MOR-T$_L$ strategy within the context of line search minimization, a commonly employed iterative process in FWI. To clarify, a line search is an iterative optimization method that minimizes a functional $J$ along a given perturbation direction, starting from an initial model $\theta_0$. In the following analysis, we consider the perturbation velocity $\dtheta_0$, which is obtained through an adjoint state method. The line search process is then applied to determine the optimal step size $\alpha^*$, such that $\theta_0 
 \alpha^* \dtheta_0$ represents an admissible updated velocity.

\begin{minipage}[t]{0.49\textwidth}
\begin{figure}[H]
    \centering
    \includegraphics[width=\textwidth]{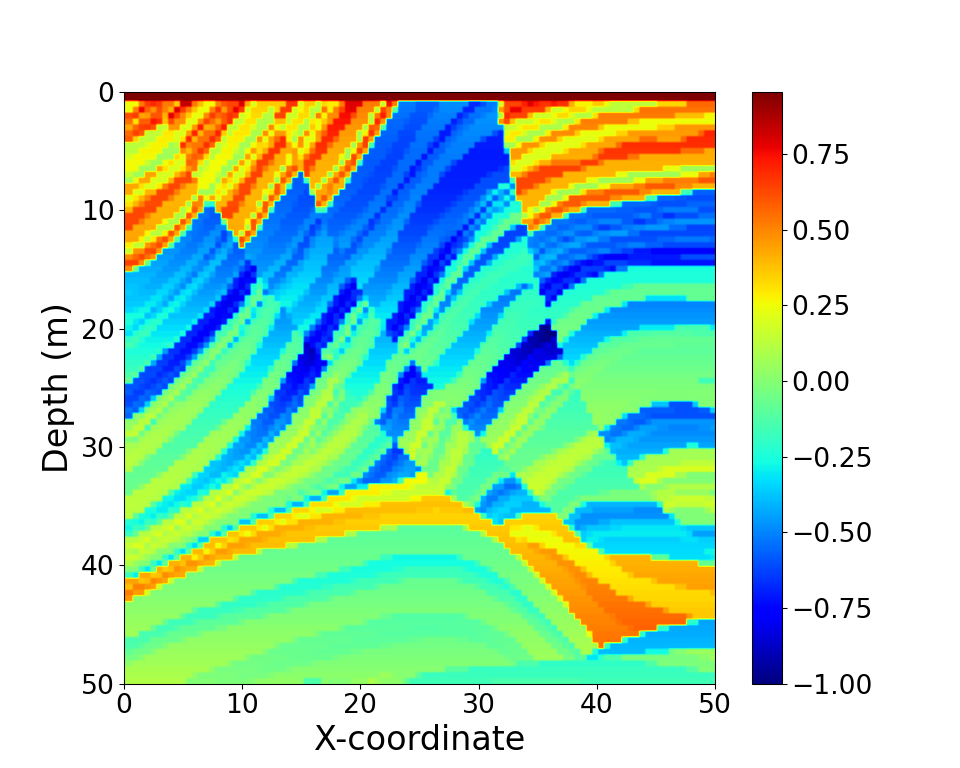} 
    \caption{{\bf{\footnotesize Normalized Perturbation Distribution.}}\\ \footnotesize \it The figure shows the normalized velocity perturbation $\dtheta_0$ in the heterogeneous medium $\Omega$, with depth represented on the y-axis and the x-coordinate on the x-axis. The color map indicates the normalized velocity perturbation values.}
    \label{fig:gradient marmousi}
\end{figure}
\end{minipage}%
\hfill
\begin{minipage}[t]{0.49\textwidth}
\begin{figure}[H]
    \centering
    \includegraphics[width=\textwidth]{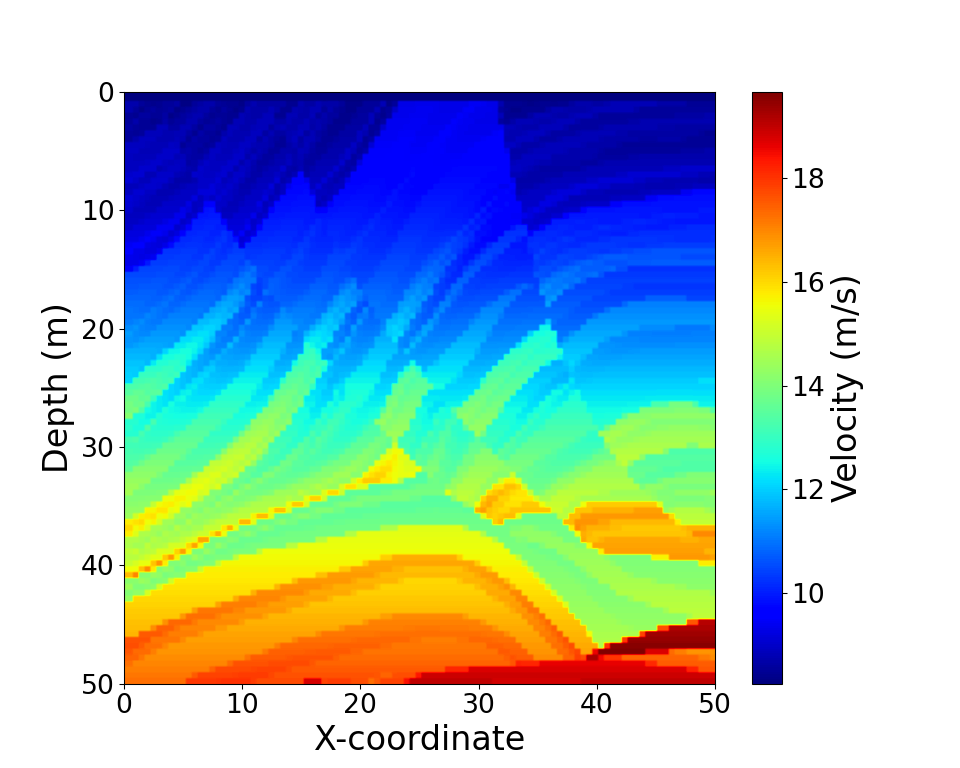} 
    \caption{{\bf{\footnotesize Updated Velocity Distribution in a Heterogeneous Medium.}} \\ \footnotesize \it The figure shows the updated velocity perturbation $c_1 = \frac{1}{\sqrt{\theta_0 + 2.12 \times 10^{-3}\dtheta_0}}$ in the heterogeneous medium $\Omega$, with depth represented on the y-axis and the x-coordinate on the x-axis. The color map indicates the velocity values.}
    \label{fig:model marmousi m1}
\end{figure}
\end{minipage}%
\bigskip

\subsubsection{Accuracy of the MOR-T$_L$ in Computing Velocity Perturbations}
\begin{figure}[H]
    %\centering
    \hspace{-0.8cm}
    \includegraphics[width=1.1\textwidth]{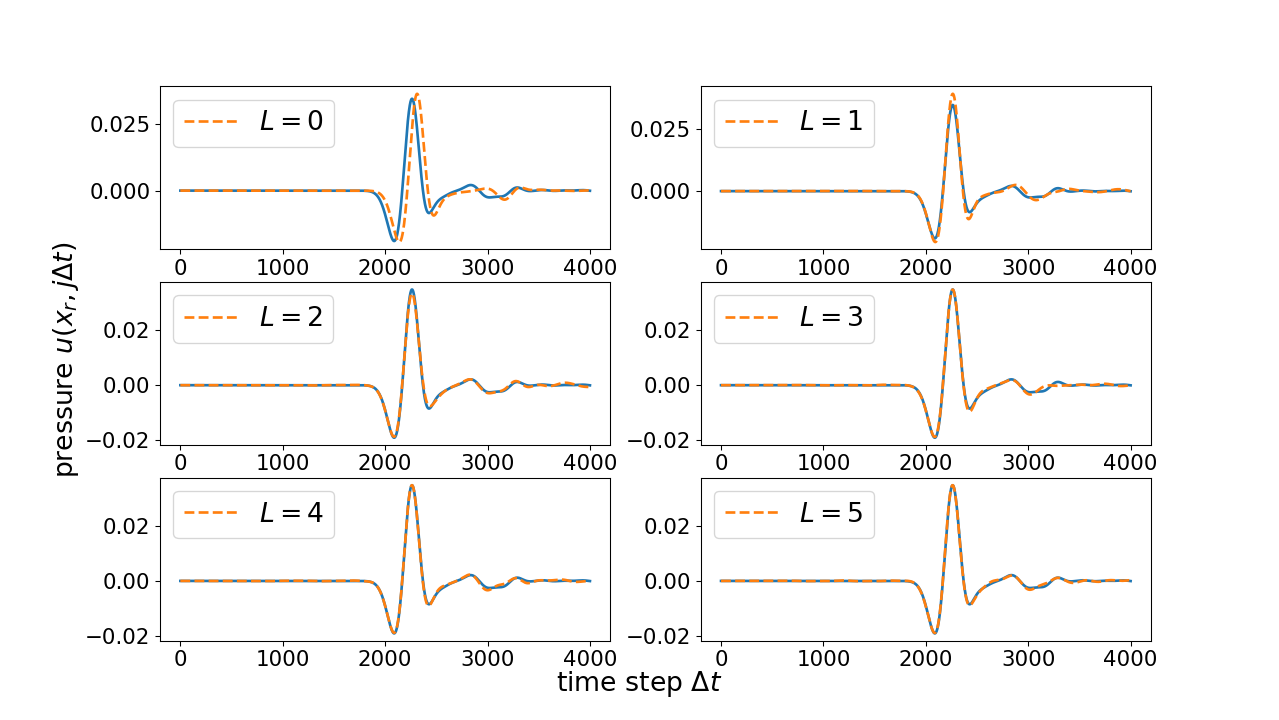} 
    \caption{{\bf{\footnotesize Accuracy of MOR-T$_L$ and Sensitivity to Taylor Polynomial Degree.}} \\ \footnotesize \it This figure consists of six subfigures illustrating the accuracy of the proposed MOR-T$_L$ method and its sensitivity to the degree $L$ of the Taylor polynomial expansion (from $L=0$ to $L=5$). Each subfigure shows a comparison of the pressure field at $x_r = (39.05, 15.42)$ over the time step interval $[0, 4\times10^3]$ for the updated velocity $c_1$ with results from SEM (reference) and MOR-T$_L$.}
    \label{fig:traces Fréchet abc marmousi m1 SEM QR}
\end{figure}

\begin{figure}[H]
    \centering
    \includegraphics[width=\textwidth]{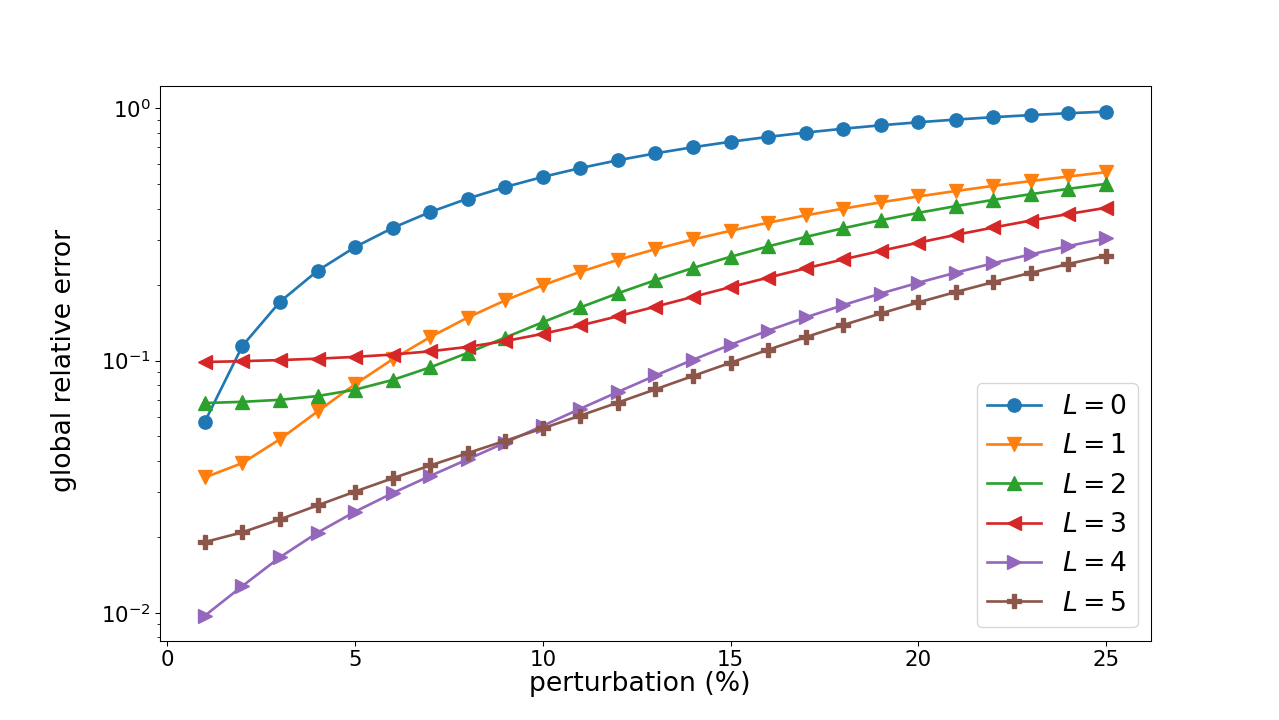}
    \caption{{\bf{\footnotesize Sensitivity of MOR-T$_L$ Accuracy to Perturbation Size for Varying Taylor Polynomial Degree L.}} \\\footnotesize \it This figure shows the sensitivity of MOR-T$_L$ accuracy to the size of the perturbation $\alpha \dtheta$. The x-axis represents the perturbation size in percentage, while the y-axis shows the average relative errors. Six curves are plotted, each corresponding to a different Taylor polynomial degree L (from L=0 to L=5). The results are computed for the updated velocity $c_1$.}
    \label{fig:error Fréchet abc relative to perturbation}
\end{figure}

To perform a line search minimization, we consider the perturbation $\dtheta_0$ as shown in Figure \ref{fig:gradient marmousi}. The updated velocity takes the form $\theta_1 = \theta_0 + \alpha \dtheta_0$, with $\alpha = 2.12 \times 10^{-3}$ resulting in a perturbation size $\alpha \frac{\lVert \dtheta_0 \rVert_2}{\lVert \theta_0 \rVert_2}$ of approximately 10\%. This perturbation represents a substantial change in the initial velocity (see Figure \ref{fig:model marmousi m1} for the updated $\theta_1$, compared with the original velocity $\theta_0$ in Figure \ref{fig:model m0 marmousi}).
We then apply the proposed MOR-T$_L$ strategy to compute the pressure field $u_{\theta_1}$ and compare the results with those from the standard SEM approach, which serves as the reference solution. In particular, we investigate the sensitivity of MOR-T$_0$’s accuracy to the degree $L$ of the Taylor polynomial (see Equation \eqref{taylor}). The results are summarized in Figures \ref{fig:traces Fréchet abc marmousi m1 SEM QR}-\ref{fig:error Fréchet abc relative to perturbation} and Table \ref{tab:relative error Fréchet}. Key observations include:

\begin{itemize}
    \item Figure \ref{fig:traces Fréchet abc marmousi m1 SEM QR} compares the pressure field at the point $\x_r = (39.05, 15.42)$ over the time interval $[0, 4\times10^3]$ for Taylor polynomial of different degrees. When $L=0$, corresponding to the standard MOR method, a phase-shift error is observed relative to the reference solution. This discrepancy is concerning, as it could significantly impact the convergence of the line search procedure. As expected, increasing the degree $L$ improves accuracy. Remarkably, the MOR-T$_L$ method achieves excellent accuracy with as few as four terms in the Taylor expansion ($L=3$).
    \item Figure \ref{fig:error Fréchet abc marmousi m1 SEM QR} shows the relative error in the L$^2$-norm over the entire spatial domain $\Omega$ and the time step interval. Unlike Figure \ref{fig:traces Fréchet abc marmousi m1 SEM QR}, which focuses on local errors, Figure \ref{fig:error Fréchet abc marmousi m1 SEM QR} presents global error estimates. The results confirm that, for $L=0$, the standard ROM method produces unacceptable errors, while the MOR-T$_L$ method achieves high accuracy even for small $L$ values, particularly for $L=3$.
    \item Figure \ref{fig:Frechet size abc marmousi m1 SEM QR} and Table \ref{tab:relative error Fréchet} provide an initial assessment of the computational cost of the proposed MOR-T$_L$ method. The number of basis functions increases linearly with $L$, and for $L=3$, the offline computational time is about 4 minutes, with a 6\% accuracy level.
    \item Figure \ref{fig:error Fréchet abc relative to perturbation} indicates that the MOR-T$_L$ method maintains acceptable accuracy even when the perturbation $\alpha \dtheta_0$ reaches 15\%, a relatively high value for the perturbation. This suggests that the proposed method is robust to large perturbations.
\end{itemize}

\begin{minipage}{0.49\textwidth}
    \centering
    \includegraphics[width=\textwidth]{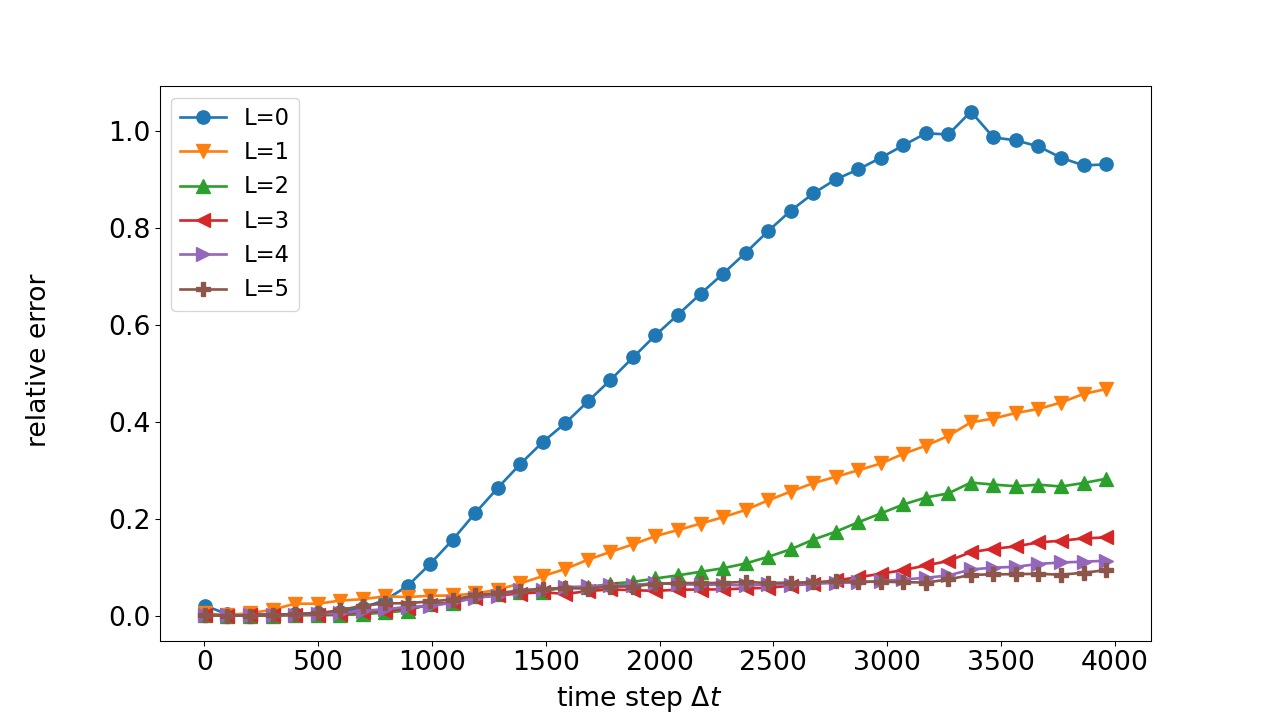}
    \caption{{\bf{\footnotesize Relative L2-norm Errors for MOR-T$_L$ with Varying Taylor Polynomial Degree.}} \\\footnotesize \it This figure presents six curves, each representing the relative L$^2$-norm error for MOR-T$_L$ with different degrees $L$ of the Taylor polynomial expansion (from $L=0$ to $L=5$). The errors are computed across the computational domain $\Omega$ at each time step in the interval $[0, 4\times10^3]$  for the updated velocity $c_1$.}
    \label{fig:error Fréchet abc marmousi m1 SEM QR}
\end{minipage}
\hfill
\begin{minipage}{0.49\textwidth}
    \centering
    \includegraphics[width=\textwidth]{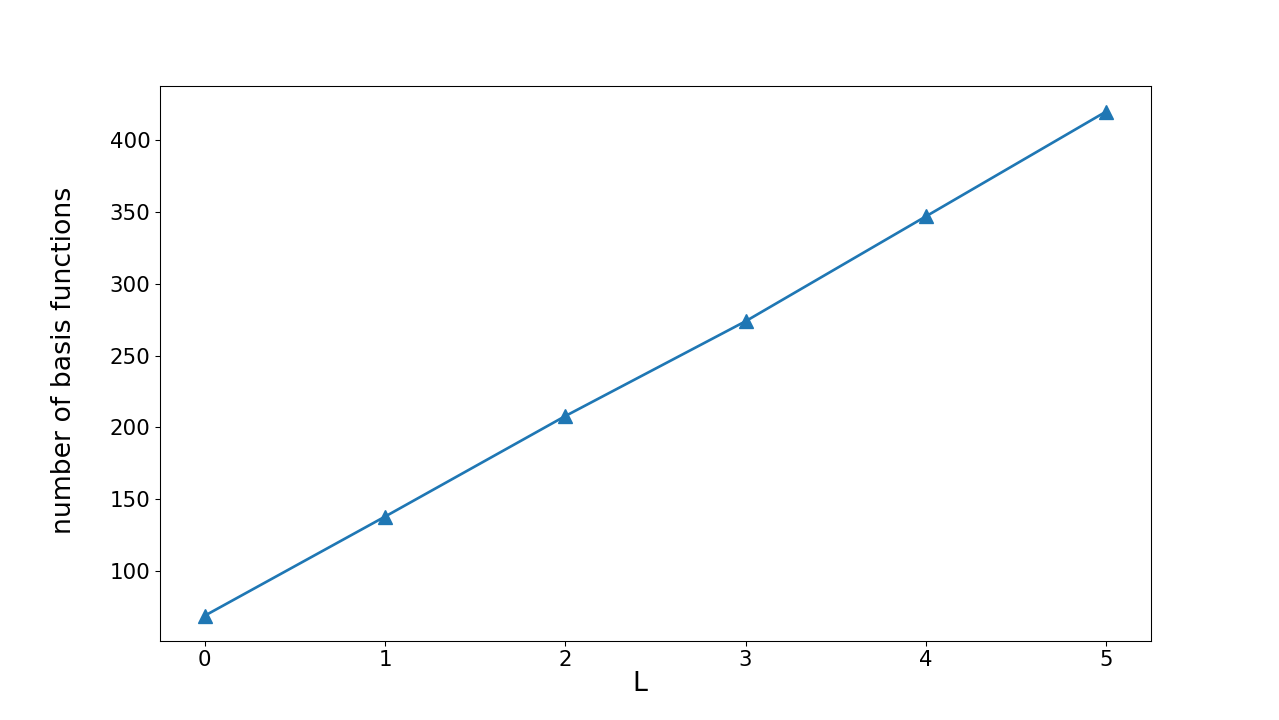}
    \caption{{\bf{\footnotesize Dependence of MOR-T$_L$ Basis Functions on Taylor Polynomial Degree L.}} \\\footnotesize \it This figure illustrates the linear relationship between the number of basis functions required by the MOR-T$_L$ method and the degree $L$ of the Taylor polynomial expansion. As $L$ increases, the number of basis functions grows linearly.}
    \vspace{0.67cm}
    \label{fig:Frechet size abc marmousi m1 SEM QR}
\end{minipage}

\bigskip

\begin{table}[H]
\scriptsize
    \centering
    \begin{tabular}{|c|c|c|c|c|c|c|}
        \hline
        $L$ &Total Computational Time  &$\displaystyle \frac{1}{N_t}\sum_{i=1}^{N_t} \frac{\lVert u_{SEM}(\x,t_i) - u_{ROM}(\x,t_i) \rVert_2}{\lVert u_{SEM}(\x,t_i) \rVert_2}$\\
        \hline
        0 &47s &0.535 \\
        \hline
        1 &120s &0.198 \\
        \hline
        2 &182s &0.113 \\
        \hline
        3 &254s &0.061 \\
        \hline
        4 &347s &0.055 \\
        \hline
        5 &561s &0.053 \\
        \hline
    \end{tabular}
    \caption{{\bf{\footnotesize MOR-T$_L$: Total Computational Time and Average Relative Errors for Varying Taylor Polynomial Degree L.}} \\\footnotesize \it This table compares the total computational time (offline time for constructing the basis and online time for solving the system) and the average relative L$^2$-norm error for MOR-T$_L$ across different values of $L$, the Taylor polynomial degree. The errors are computed over the computational domain $\Omega$ and averaged over time steps for the updated velocity $c_1$.}
    \label{tab:relative error Fréchet}
\end{table}

In summary, these results emphasize the inadequacy of the standard ROM approach for accurately computing pressure fields at updated velocities. The proposed MOR-T$_L$ method, however, delivers high accuracy even with a cubic Taylor expansion ($L=3$), and its computational cost remains manageable due to the linear growth of the basis functions with $L$.

\subsubsection{Computational Efficiency and Robustness of the MOR-T$_L$ During Line Search}

In this section, we evaluate the computational efficiency and robustness of the MOR-T$_L$ method in the line search, particularly in the presence of noisy measurements. Specifically, the measurements of the pressure field vector $u_{obs}$ is obtained by computing the pressure $u_{\theta^*}$ corresponding to the updated model parameter $\theta_0 + \alpha^* \dtheta_0$ with $\alpha^*= 2.12 \times 10^{-3}$ over the time interval $[0, 4\times10^3]$ on a line of $202$ equidistant receivers from $\x_{r_1}=(0,15.42)$ to $\x_{r_{202}}=(50,15.42)$. 

These synthetic measurements are obtained using the SEM method to avoid the so-called inverse crime. We then compute the pressure fields $u_{\theta_\alpha}$ for a set of velocities $\theta_\alpha = \theta_0 + \alpha\dtheta_0$ using the MOR-T$_L$ method.
The goal of the line search is to minimize the cost function $J(\alpha)$ defined as: 

\begin{equation}\label{cost function}
    \begin{aligned}
        J(\alpha) = \frac{1}{2}\lVert u_{\theta_\alpha} - u_{obs} \rVert^2_2 
    \end{aligned}
\end{equation}
by evaluating it at several values of $\alpha$. We also introduce white noise into the vector measurements $u_{obs}$ (see Figure \ref{fig:cost functions}) to assess the method’s robustness under noisy conditions. The results are presented in Figures \ref{fig:cost functions} and Tables \ref{tab:relative error Fréchet}-\ref{tab:relative computational Time Frechet}, with the following key findings:

\begin{itemize}
    \item Preliminary Remark: The line search is conducted by evaluating $J(\alpha)$ at 200 points within the interval $[0, 4\times10^{-3}]$ which is uniformly subdivided. This interval is chosen to account for perturbations $\alpha\dtheta_0$ of up to 20\% on the relative errors, representing relatively large perturbations. Negative values of $\alpha$ are excluded to ensure that the perturbation $\alpha\dtheta_0$ aligns with the direction of $-\nabla_\theta{J}$, as required by the steepest descent method.
    \item Figure \ref{fig:cost functions}(a) shows results for noise-free measurements. The standard MOR method ($L=0$) fails to identify the optimal $\alpha^*$, as expected given the sensitivity of the ROM basis to model parameters. In contrast, the MOR-T$_L$ method accurately determines $\alpha^*$ even when using a linear Taylor expansion ($L=1$). For all non-constant Taylor polynomials, $J(\alpha)$ exhibits a distinct minimum around the correct value of $\alpha^*$, facilitating an unambiguous solution.
    \item Figure \ref{fig:cost functions}(b)-(c) show results for noise levels of $1\%$ and $5\%$, respectively. Despite the noise, the MOR-T$_L$ method remains robust and continues to determine $\alpha^*$ with impressive accuracy, achieving a relative error of less than $1\%$ in both cases, as reported in Table \ref{tab:optimal alpha}.
    \item Table \ref{tab:relative computational Time Frechet} compares the computational costs of the MOR-T$_L$ and SEM methods when evaluating $J(\alpha)$ at $10$ parameter values, consistent with the use of Wolfe conditions in FWI to limit the number of evaluations. We observe that MOR-T$_3$ reduces the computational time by about 50\%. While evaluating $J(\alpha)$ at a small number of points may compromise the accuracy of $\alpha^*$ and slow down convergence, the MOR-T$_L$ method still outperforms SEM in terms of computational cost and accuracy.
    \item Final Remark: The use of SEM in the line search stage can result in prohibitive computational costs when a large number of parameter points are needed. For example, evaluating 200 points with SEM requires over three hours of computation, while MOR-T$_L$ with a Taylor polynomial degree $L$ smaller than $5$ completes the task in under 10 minutes. This efficiency allows for the use of many points during the line search, which is crucial for ensuring fast convergence and accuracy in the FWI process. In contrast, relying on only a few points based on the Wolfe conditions -- especially if the initial guess is outside the pre-asymptotic region -- can significantly hinder convergence, or even prevent it entirely. Consequently, any computational savings from using fewer points may be offset by the need for more iterations, ultimately negating the intended efficiency gains.
\end{itemize}

\begin{minipage}[t]{\textwidth}
\hfill
\begin{minipage}[t]{0.49\linewidth}
    \centering
    \includegraphics[width=\linewidth]{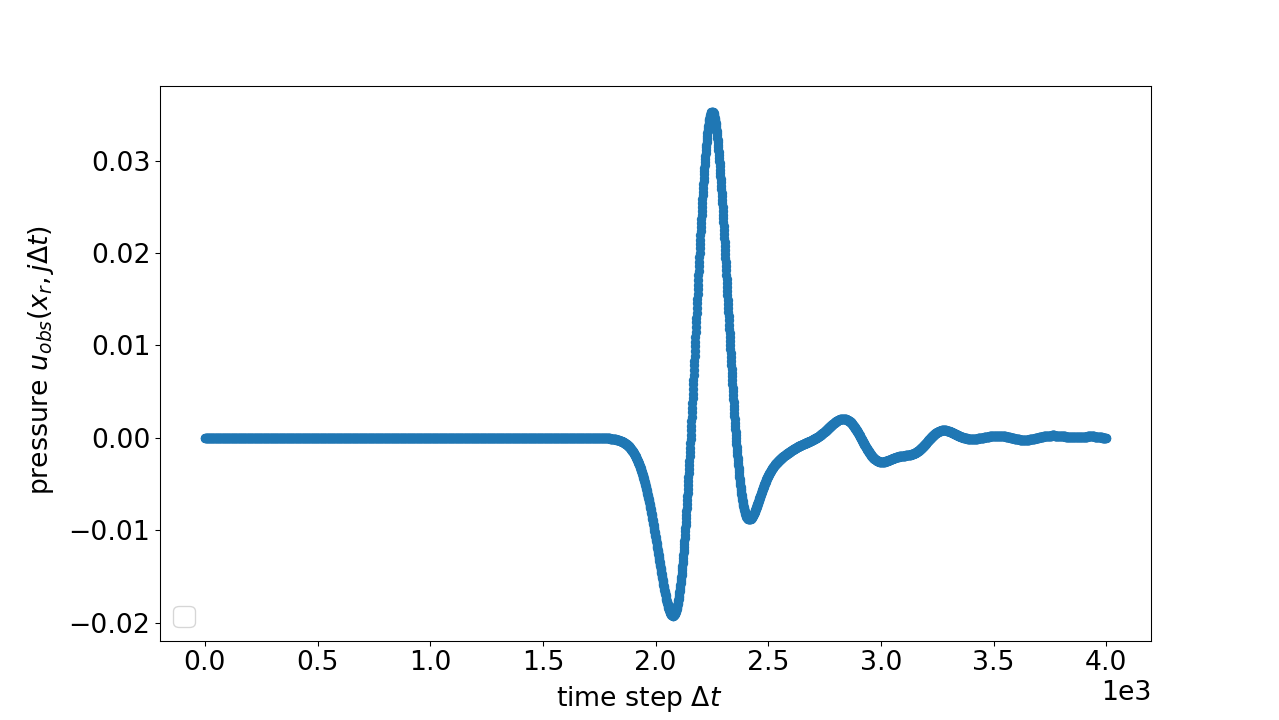}
\end{minipage}%
\hfill
\begin{minipage}[t]{0.49\linewidth}
    \centering
    \includegraphics[width=\linewidth]{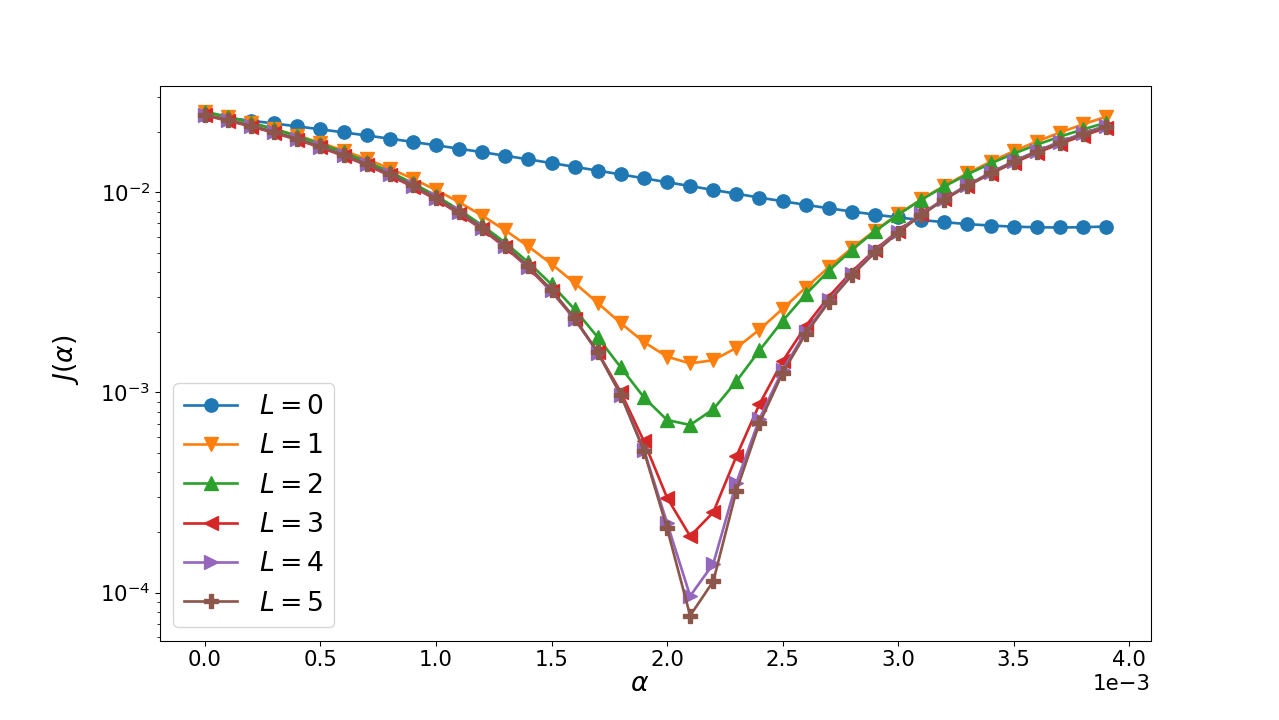}
\end{minipage}
\centering{(a) \footnotesize Measurement (left) and misfit function $J(\alpha)$ (right) for noise level: $0\%$}
\end{minipage}
\vspace{1em}
\begin{minipage}[t]{\textwidth}
\hfill
\begin{minipage}[t]{0.49\textwidth}
    \centering
    \includegraphics[width=\textwidth]{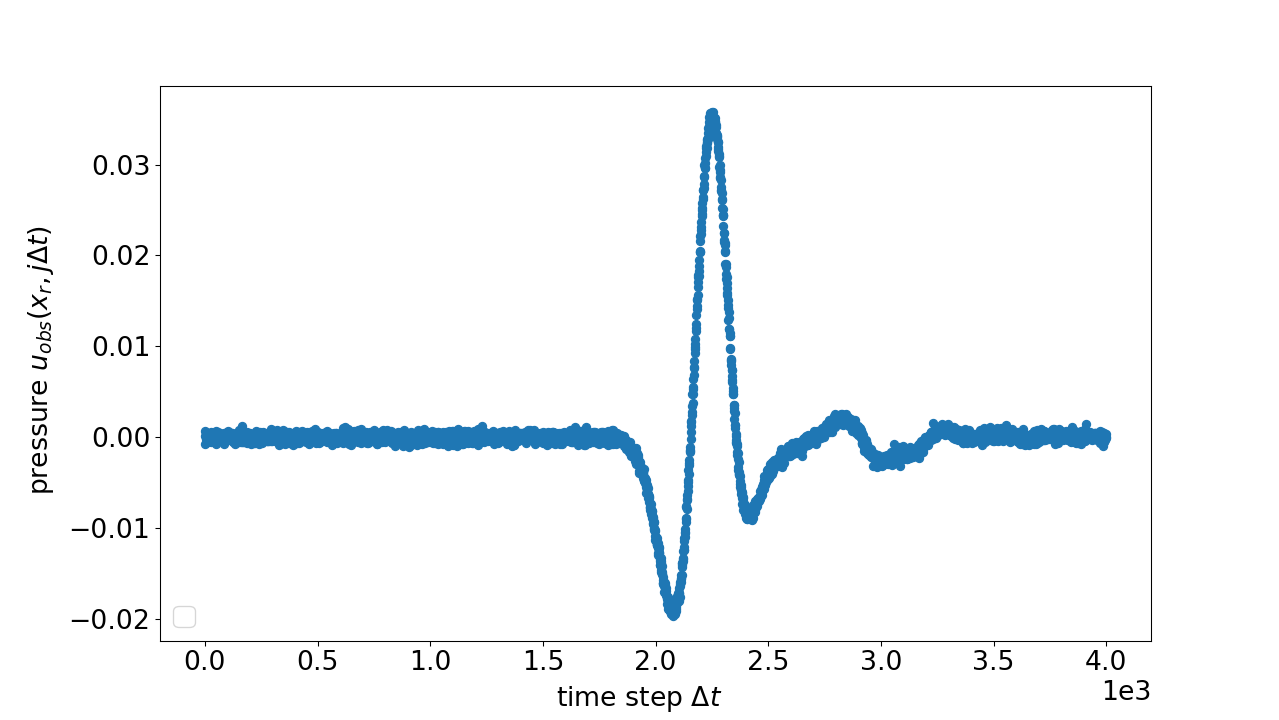}
\end{minipage}%
\hfill
\begin{minipage}[t]{0.49\textwidth}
    \centering
    \includegraphics[width=\textwidth]{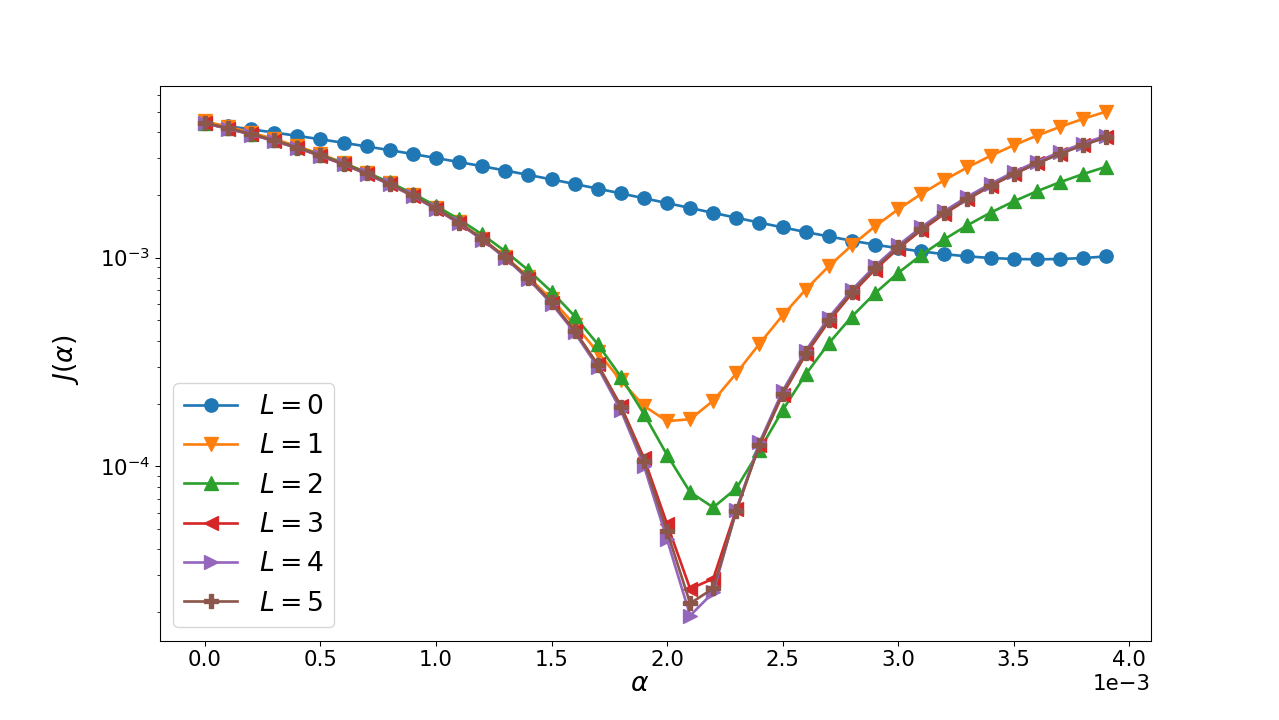}
\end{minipage}
\centering{(b) \footnotesize Measurement (left) and misfit function $J(\alpha)$ (right) for noise level: $1\%$}
\end{minipage}
\vspace{1em}
\begin{minipage}[t]{\textwidth}
\hfill
\begin{minipage}[t]{0.49\textwidth}
    \centering
    \includegraphics[width=\textwidth]{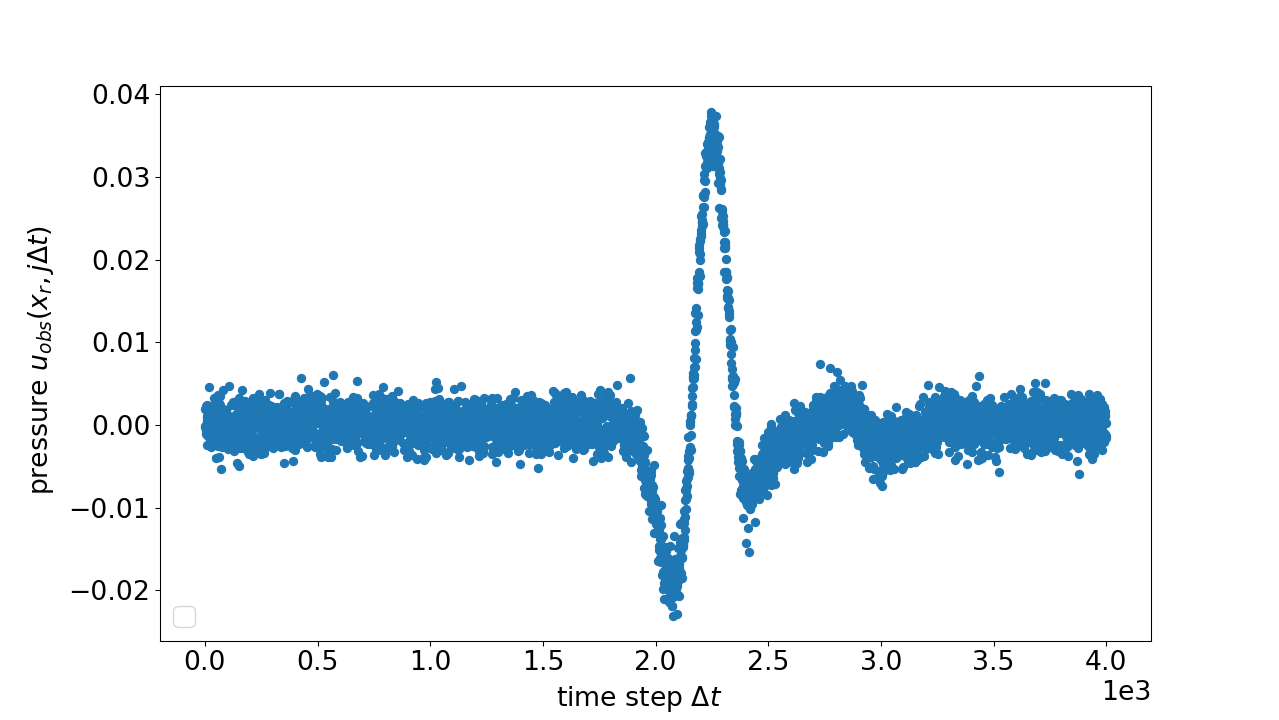}
\end{minipage}%
\hfill
\begin{minipage}[t]{0.49\textwidth}
    \centering
    \includegraphics[width=\textwidth]{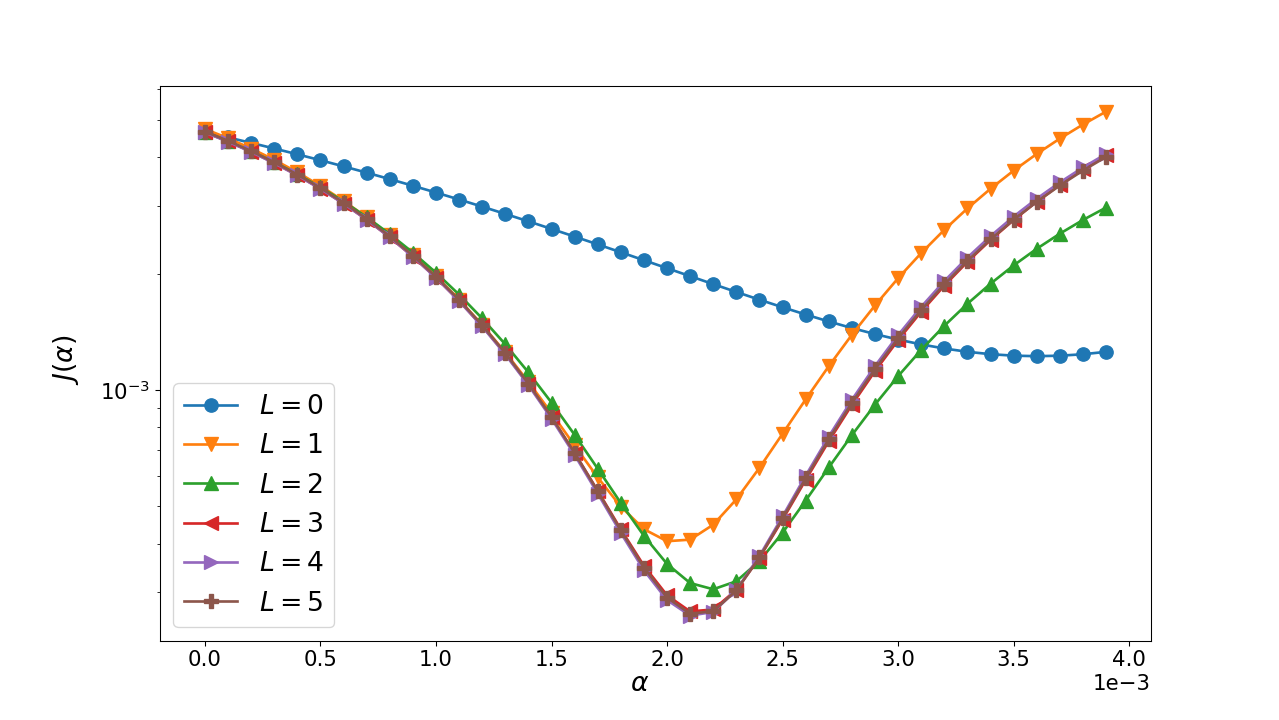}
\end{minipage}
\centering{(c) \footnotesize Measurement (left) and misfit function $J(\alpha)$ (right) for noise level: $5\%$}
\end{minipage}
\captionof{figure}{{\bf{\footnotesize Pressure Field Measurements and Cost Function with Varying Noise Levels.}} \\ \footnotesize \it This figure consists of six subfigures. The left three show pressure field measurements at $x_r = (39.05, 15.42)$ over $[0, 4\times10^3]$ for the velocity $c^*$ obtained with SEM and corrupted with 0\%, 1\%, and 5\% white noise (top to bottom). The right three subfigures display the cost function $J(\alpha)$ evaluated at 200 parameter values, with six curves corresponding to Taylor polynomial degrees $L=0$ to $L=5$ for each noise level.}
\label{fig:cost functions}

\bigskip

\begin{table}[H]
    \centering
    \begin{tabular}{|c|c|c|c|}
    \hline
        & \multicolumn{3}{|c|}{$argmin$ $J(\alpha)$ for noise level} \\
        \hline
        $L$ & 0\% &1\% &5\% \\
        \hline
        0 &$3.7 \times 10^{-3}$ &$3.6\times 10^{-3}$ &$3.6\times 10^{-3}$ \\
        \hline
        1 &$2.12\times 10^{-3}$ &$2.04\times 10^{-3}$ &$2.04\times 10^{-3}$ \\
        \hline
        2 &$2.08\times 10^{-3}$ &$2.2\times 10^{-3}$ &$2.2\times 10^{-3}$ \\
        \hline
        3 &$2.12\times 10^{-3}$ &$2.14\times 10^{-3}$ &$2.14\times 10^{-3}$ \\
        \hline
        4 &$2.12\times 10^{-3}$ &$2.14\times 10^{-3}$ &$2.14\times 10^{-3}$ \\ 
        \hline
        5 &$2.12\times 10^{-3}$ &$2.14\times 10^{-3}$ &$2.14\times 10^{-3}$ \\
        \hline
    \end{tabular}
    \caption{{\bf{\footnotesize Retrieved Parameter $\alpha^*$ by MOR-T$_L$ for Varying Noise Levels and Taylor Polynomial Degree L.}} \\\footnotesize \it This table shows the computed values of $\alpha^*$ for different noise levels (0\%, 1\%, and 5\%) and Taylor polynomial degrees $L$. It illustrates the effect of noise and $L$ on the accuracy of $\alpha^*$ and the robustness of the MOR-T$_L$ strategy.}
    \label{tab:optimal alpha}
\end{table}

\begin{table}[H]
    \centering
    \begin{tabular}{|c|c|}
        \hline
        $L$ &$\frac{\text{MOR-T}_L \text{ Total Computational Time}}{\text{SEM Computational Time}}$\\
        \hline
        0  &\centering 0.12 \tabularnewline
        \hline
        1  &\centering 0.27 \tabularnewline
        \hline
        2  &\centering 0.39 \tabularnewline
        \hline
        3  &\centering 0.55 \tabularnewline
        \hline
        4  &\centering 0.74 \tabularnewline
        \hline
        5  &\centering 1.16 \tabularnewline
        \hline
    \end{tabular}
    \caption{{\bf{\footnotesize Computational Time Ratio between MOR-T$_L$ and SEM for Line Search with 10 Parameter Values.}} \\\footnotesize \it This table shows the time ratio between MOR-T$_L$ and SEM during a line search with 10 parameter values $\alpha$, highlighting the computational cost efficiency of MOR-T$_L$ for different Taylor polynomial degrees L.}
    \label{tab:relative computational Time Frechet}
\end{table}

In conclusion, the results confirm that the standard MOR method is inadequate for line search due to its poor performance in determining the optimal $\alpha^*$. The MOR-T$_L$ method, however, excels in both accuracy and computational efficiency, making it a strong candidate for practical FWI applications. Furthermore, the MOR-T$_L$ method remains robust even in the presence of noise, offering significant advantages over SEM, particularly in terms of computational feasibility. 

\section{Summary and Conclusion}
In this work, we introduced a novel MOR approach, called MOR-T$_L$, based on Taylor polynomial expansions. The key innovation of MOR-T$_L$ lies in the construction of the ROM basis functions, which are generated using Taylor polynomial representations. This approach provides a more robust method for handling parameter variations effectively. Numerical experiments on a two-dimensional wave prototype problem, relevant to seismic applications, demonstrated the robustness of the MOR-T$_L$ basis to parameter changes and noise. Specifically, MOR-T$_L$ excels in both accuracy and computational efficiency, particularly in a line search minimization context, and can accommodate parameter perturbations of up to 20\%, which in the presence of noise, offers substantial advantages over SEM in terms of computational feasibility. These results position MOR-T$_L$ as a promising candidate for practical applications in FWI, demonstrating its potential to improve computational cost-effectiveness and robustness.

% ---------------------------------------------------------------------

\appendix

%%%%%%%%%% Bibliography
%% \newpage 
%\footnotesize 
\bibliographystyle{siam}
\bibliography{article}

\end{document}